\documentclass{raex}

\begin{Author} 
	\FirstName{Tam\'as}\LastName{M\'atrai}
	\PostalAddress{University of Toronto, 40 St George St., Toronto, ON, Canada, M5S 2E4}
	\Email{matrait@renyi.hu}
	\Thanks{The research was partially supported  by the OTKA grants K 61600, K 49786 and K 72655. We thank M\'at\'e Vizer for the inspiring discussions. We also thank Longyun Ding and the referee for pointing out an error in an earlier version of this paper.}
\end{Author}

\begin{MathReviews} 
	\primary{03E15}
	\secondary{46A45}
\end{MathReviews}

\begin{KeyWords} 
\end{KeyWords}

\title{On $\ell^{p}$-like equivalence relations}

\usepackage{amsfonts,amssymb,amsmath}
\usepackage[mathscr]{eucal}
\usepackage{graphicx,enumerate}

\begin{document} 
 
\maketitle   

\def\es{\emptyset} 
\def\ve{\varepsilon} 
\def\vp{\varphi} 
\def\vt{\vartheta} 
\def\sm{\setminus} 
\def\real{\mathbb{R}}
\def\CG{\mathbb{S}} 
\def\cmpl{\mathbb{C}} 
\def\cl{\mathrm{cl}} 
\def\bel{\mathrm{int}} 
\def\egesz{\mathbb{Z}} 
\def\nat{\mathbb{N}} 
\def\b{\left} 
\def\j{\right} 
\def\<{\left<} 
\def\>{\right>} 
\def\rar{\rightarrow} 
\def\lrar{\leftrightarrow} 
\def\ben{\begin{enumerate}} 
\def\een{\end{enumerate}} 
\def\Keq{\begin{equation}} 
\def\Zeq{\end{equation}} 
\def\Kem{\begin{multline}} 
\def\Zem{\end{multline}} 
\def\bit{\begin{itemize}} 
\def\eit{\end{itemize}} 
\def\pa{\partial} 
\def\d{\mathrm{d}} 
\def\cpl{\complement} 
\def\cf{\mathrm{cf}} 
\def\bE{\mathbf{E}} 
\def\Id{\mathrm{Id}}
\def\card{\mathrm{card}} 
\def\ss{\subseteq} 
\def\Dh{\Delta_{h}} 
\def\vv{\vert \vert} 
\def\mc{\mathcal} 
\def\DS{\displaystyle} 
\def\TS{\mathstyle} 
\def\BS{{\bf \Sigma}} 
\def\BP{{\bf \Pi}} 
\def\sign{\mathrm{sign}} 
\def\diam{\mathrm{diam}} 
\def\dist{\mathrm{dist}} 
\def\PR{\textstyle \mathrm{Pr}} 
\def\bs{\blacksquare} 
\def\prf{\textbf{Proof. }}

\begin{abstract} For $f \colon [0,1] \rar \real^{+}$, consider the relation $\mathbf{E}_{f}$ on $[0,1]^{\omega}$ defined by $(x_{n}) \mathbf{E}_{f} (y_{n}) \Leftrightarrow \sum_{n < \omega} f(|y_{n} - x_{n}|) < \infty.$ We study the Borel reducibility of Borel equivalence relations of the form $\mathbf{E}_{f}$. Our results indicate that for every $1 \leq p < q < \infty$, the order $\leq_{B}$ of Borel reducibility on the set of equivalence relations $\{\bE \colon \bE_{\Id^{p}} \leq_{B} \bE \leq _{B} \bE_{\Id^{q}}\}$ is more complicated than expected, e.g.\ consistently every linear order of cardinality continuum embeds into it.
\end{abstract}

\newtheorem{theorem}{Theorem}   
\newtheorem{corollary}[theorem]{Corollary} 
\newtheorem{lemma}[theorem]{Lemma} 
\newtheorem{clm}[theorem]{Proposition} 
\newtheorem{definition}[theorem]{Definition} 
\newtheorem{que}[theorem]{Question} 
\newtheorem{megj}[theorem]{Remark} 
\newtheorem{prob}[theorem]{Problem}

\section{Introduction} 

Let $f \colon [0,1] \rar \real^{+}$ be an arbitrary function and consider the relation $\mathbf{E}_{f}$ on $[0,1]^{\omega}$ defined by setting, for every $(x_{n})_{n < \omega}, (y_{n})_{n < \omega} \in [0,1]^{\omega}$, \Keq\label{def} (x_{n}) \mathbf{E}_{f} (y_{n}) \Leftrightarrow \sum_{n < \omega} f(|y_{n} - x_{n}|) < \infty.\Zeq Several natural questions arise, e.g.\
\bit
\item[(i)] when is $\mathbf{E}_{f}$ an equivalence relation?
\item[(ii)] which equivalence relations can be obtained in the form $\mathbf{E}_{f}$?
\item[(iii)] for what $f,g \colon [0,1] \rar [0,1]$ is $\mathbf{E}_{f}$ Borel reducible to $\mathbf{E}_{g}$?
\eit In the present paper we answer (i), we initiate a study of (ii) and we obtain various conditions for (iii).
 
The prototypes of equivalence relations of the form $\mathbf{E}_{f}$ are induced by the Banach spaces $\ell^{p}$ $(1 \leq p < \infty)$, i.e.\ they are defined by the functions $f = \mathrm{Id}^{p}$ for $1 \leq p < \infty$, where $\mathrm{Id} \colon [0,1] \rar [0,1]$ is the identity function. The Borel reducibility among these equivalence relations is fully described by a classical result of R.\ Dougherty and G.\ Hjorth \cite[Theorem 1.1 p.\ 1836 and Theorem 2.2 p.\ 1840]{DH} stating that for every $1 \leq p, q < \infty$,  \Keq\label{dh}\mathbf{E}_{\mathrm{Id}^{p}} \leq _{B} \mathbf{E}_{\mathrm{Id}^{q}}  \Leftrightarrow p \leq q.\Zeq We note, however, that e.g.\ for the function $f(0) = 0$, $f(x) = 1$ $(0 < x \leq 1)$ we have $\mathbf{E}_{f}$ is the equivalence relation of eventual equality on $[0,1]^{\omega}$, also denoted by $E_{1}$ in the literature; that is, the investigation of equivalence relations of the form $\mathbf{E}_{f}$ concerns equivalence relations which are not necessarily reducible to $\bE_{\Id^{p}}$ for some $1 \leq  p < \infty$.
 
Our investigations were motivated by a question of S.\ Gao in \cite{G} p.\ 74, asking whether for $1 \leq p < \infty$, $\mathbf{E}_{\mathrm{Id}^{p}}$ is the greatest lower bound of $\{\mathbf{E}_{\mathrm{Id}^{q}} \colon p < q < \infty\}$; we note that formally the question in \cite{G} p.\ 74 refers to equivalence relations on $\real^{\omega}$, but as we will see later in Lemma \ref{mege}, the two formulations are equivalent.
We answer this question in the negative by showing, for fixed $1 \leq p  < \infty$, that $\mathbf{E}_{\mathrm{Id}^{p}} < _{B} \bE_{f} <_{B} \mathbf{E}_{\mathrm{Id}^{q}}$ for every $q > p$ whenever  \Keq\label{nov}\lim_{x \rar +0}\frac{f(x)}{x^{p}} = 0 \textrm{ and } \lim_{x \rar +0}\frac{f(x)}{x^{q}} = \infty  ~(p < q < \infty),\Zeq and $f$ satisfies some additional technical assumptions (see e.g.\ Corollary \ref{megm}). However, toward this result we aim to carry out a general study of the relations $\mathbf{E}_{f}$ and their Borel reducibility. To this end, in Section \ref{gen} we characterize the functions for which $\bE_{f}$ is an equivalence relation and, roughly speaking, we show that $f$ is continuous if and only if $E_{1} \not \leq _{B} \bE_{f}$. In Section \ref{red} and in Section \ref{nonred} we prove general reducibility and nonreducibility results for equivalence relations of the form $\mathbf{E}_{f}$. The results of these sections heavily build on techniques developed in \cite{DH}. Finally, in Section \ref{con} we conclude our investigations by applying the technical results of the previous sections to concrete functions; in particular, we answer the above mentioned question of S.\ Gao, and we show that for $1 \leq p < q < \infty$, every linear order which embeds into $(\mc{P}(\omega)/\mathrm{fin},\subset)$ also embeds into the set of equivalence relations $\{\bE_{f} \colon \bE_{\Id^{p}} \leq_{B} \bE_{f} \leq _{B} \bE_{\Id^{q}}\}$ ordered by $<_{B}$. 

Our results produce just examples. We are far from giving a full description of the Borel equivalence relations of the form $\bE_{f}$ or a complete picture of the Borel reducibility relation among the $\bE_{f}$s. In particular, it remains open whether  there are two functions $f$ and $g$ such that $\bE_{f}$ and $\bE_{g}$ are incomparable under $\leq_{B}$. Nevertheless, we have one qualitative observation. Conditions (\ref{dh}) and (\ref{nov}) may suggest that reducibility among the $\bE_{f}$s is essentially governed by the growth order of the $f$s. However, this is far from being true. As we will see in Section \ref{red}, under mild additional assumptions on $f$ we have e.g.\ $ \bE_{  \mathrm{Id}^{p}} \leq _{B} \bE_{f}$ whenever $\lim_{x \rar +0} f(x)/x^{p-\ve} = 0$ for every $\ve > 0$ (for the precise statement, see Theorem \ref{ala}); this is in contrast with (\ref{nov}). 

For basic terminology in descriptive set theory we refer to \cite{K}. As above, if $X$ and $Y$ are Polish spaces, $E$ and $F$ are  equivalence relations on $X$ and $Y$, then we say  \emph{$E$ is Borel reducible to $F$}, $E \leq _{B} F$ in notation, if  there exists a Borel function $\vt \colon X \rar Y$ satisfying $$x E x'   \Leftrightarrow \vt(x) F \vt(x').$$ We say  \emph{$E$ and $F$ are Borel equivalent} if $E \leq _{B} F$ and $F \leq _{B} E$, while we write $E < _{B} F$ if $E \leq _{B} F$ but $F \not \leq _{B} E$.

Depending on the context, $| \cdot |$ denotes the absolute value of a real number, the length of a sequence or the cardinality of a set; $\lfloor \cdot \rfloor$ and $\{\cdot\}$ stand for lower integer part and fractional part. We denote by $\egesz$ and $\real^{+}$ the set of integers and nonnegative reals.

\section{Basic properties} \label{gen}

\begin{definition}\rm\label{Csop}Let $(G,+)$ be an Abelian group and let $H \ss G$ satisfy 
\bit
\item[($H_{1}$)] $0 \in H$;
\item[($H_{2}$)] for every $x,y \in H$, $x-y \in H$ or $y-x \in H$; 
\item[($H_{3}$)] for every $x,y,z \in H $, $x-y \in H$ and  $y-z \in H$ implies $x-z \in H$. 
\eit
For every $x \in H \cup -H$, let $x^{+} = x$ if $x \in H$ and $x^{+} = -x$ if $x \in -H \sm H$.

 For every function $f \colon H \rar \real^{+}$, we define the relation $\mathbf{E}_{f}$ on $H^{\omega}$ by setting, for every $(x_{n})_{n < \omega}, (y_{n})_{n < \omega} \in H^{\omega}$, \Keq\label{defG} (x_{n}) \mathbf{E}_{f} (y_{n}) \Leftrightarrow \sum_{n < \omega} f((y_{n} - x_{n})^{+}) < \infty;\Zeq
the definition is valid by $(H_{2})$. 

We say $f \colon H \rar \real^{+}$ is \emph{even} if for every $x \in H \cap -H$, $f(x) = f(-x)$. 
\end{definition}

Observe that for $\tilde f \colon H \rar [0,1]$, $\tilde f(x) = \min\{f(x), 1\}$ $(x \in H)$ we have $\mathbf{E}_{\tilde f} = \mathbf{E}_{f}$. So in the sequel we only consider bounded functions.

We start this section by characterizing the bounded functions $ f \colon H \rar \real^{+}$ for which $\bE_{f}$ is an equivalence relation. To avoid a meticulous bookkeeping of non-relevant constants, we will use the terminology ``by ($\star$), $A \lesssim B$" to abbreviate that ``by property $(\star)$, there is a constant $C>0$ depending on the parameters of  $(\star)$ such that $A \leq CB$". The relations $\gtrsim$ and $\approx$ are defined analogously.

\begin{clm}\label{ekv}  Let $f \colon H \rar \real^{+}$ be a bounded even function. Let $\bE_{f}$ be the relation on $H^{\omega}$ defined by (\ref{defG}). Then $\bE_{f}$ is an equivalence relation if and only if the following conditions hold:
\bit
\item[($R_{1}$)] $f(0)=0$;
\item[($R_{2}$)] there is a $C \geq 1$ such that for every $x,y \in H$ with $x+y \in H$, 
$$\begin{array}{ll} (a) & f(x+y) \leq C(f(x)+f(y)),  \\  & \\ (b) & f(x) \leq C(f(x+y) + f(y)).\end{array}$$
\eit
\end{clm}
\prf Since $f$ is even, $\bE_{f}$ is symmetric. It is obvious that ($R_{1}$) is equivalent to $\bE_{f}$ being reflexive, so it remains to show that ($R_{2}$) is equivalent to transitivity.

Suppose first ($R_{2}$) holds and let $(x_{n})_{n < \omega}$, $(y_{n})_{n < \omega}$, $(z_{n})_{n < \omega} \in H^{\omega}$ such that $(x_{n}) \bE_{f}(y_{n})$ and $(y_{n}) \bE_{f}(z_{n})$. Let $n < \omega$ be fixed. Since the role of $x_{n}$ and $z_{n}$ is symmetric, by $(H_{2})$ we can assume $z_{n} - x_{n} \in H$.  We distinguish several cases.

If $x_{n}-y_{n} \in H$ then by $(H_{3})$, $z_{n}-y_{n} \in H$ so by $(R_{2} b)$ using $(z_{n} - y_{n}) = (z_{n} - x_{n}) + (x_{n} - y_{n})$, $$  f(z_{n} - x_{n})  \lesssim f(z_{n} - y_{n})+f((y_{n} - x_{n})^{+}).$$
If $y_{n}-x_{n} \in H$ then either $z_{n}-y_{n} \in H$, hence by $(R_{2} a)$, using $(z_{n} - x_{n}) = (z_{n}-y_{n} ) + (y_{n}-x_{n} )$, $$f (z_{n} - x_{n})  \lesssim f(z_{n} - y_{n})+f(y_{n} - x_{n});$$ or $y_{n}-z_{n} \in H$ hence by $(R_{2} b)$ using $(y_{n} - z_{n}) + (z_{n}-x_{n} ) = (y_{n}-x_{n} )$, $$f (z_{n} - x_{n})  \lesssim f((z_{n} - y_{n})^{+})+f(y_{n} - x_{n}).$$

Thus $$\sum_{n < \omega} f((z_{n} - x_{n})^{+})   \lesssim \sum_{n < \omega} f((z_{n} - y_{n})^{+}) + \sum_{n < \omega} f((y_{n} - x_{n})^{+}) < \infty,$$ which gives $(x_{n}) \bE_{f}(z_{n})$;  i.e.\ $(R_{2})$ implies transitivity.

To see the other direction, suppose first there is no $C\geq 1$ for which $(R_{2} a)$ holds, i.e.\ for every $n < \omega$ there are $\xi_{n}, \eta_{n} \in H$ such that $\xi_{n} + \eta_{n} \in H$ and \Keq \notag f(\xi_{n}+\eta_{n}) > 2^{n}  (f(\xi_{n})+f(\eta_{n})).\Zeq Set $k_{n} = \max\{1, \lfloor 1/f(\xi_{n}+\eta_{n}) \rfloor\}$;  if $B \geq 1$ is an upper bound of $f$, we have \Keq\label{1}B \geq k_{n}f(\xi_{n}+\eta_{n})\geq \frac{1}{2} \textrm{ and } B > 2^{n} k_{n} (f(\xi_{n})+f(\eta_{n})).\Zeq Let $ (x_{m})_{m < \omega} \in H^{\omega}$ be the sequence which, for every $n < \omega$, admits the value $\xi_{n}$ with multiplicity $k_{n}$; and define the sequence $ (y_{m})_{m < \omega} \in H^{\omega}$ to admit $\eta_{n}$ exactly there where $ (x_{m})_{m < \omega}$ admits $\xi_{n}$ $(n < \omega)$. Then by (\ref{1}), $$\sum_{m < \omega} f(x_{m}) < 2B \textrm{ and } \sum_{m < \omega} f(y_{m}) < 2B,$$ i.e. if $\underline{0}$ denotes the constant zero sequence we have $\underline{0} \bE_{f} (x_{m})$ and $(x_{m}) \bE_{f} (x_{m}+y_{m})$. Also by (\ref{1}), $$\sum _{m < \omega} f(x_{m} + y_{m}) = \infty,$$ i.e.\ $\underline{0} \not \!\!\bE_{f} (x_{m}+y_{m})$, which shows transitivity fails. 

Finally suppose there is no $C\geq 1$ for which $(R_{2} b)$ holds, i.e.\ for every $n < \omega$ there are $\xi_{n}, \eta_{n} \in H$ such that \Keq \notag f(\xi_{n}) > 2^{n}  (f(\xi_{n}+\eta_{n})+f(\eta_{n})).\Zeq Set $k_{n} = \max\{1, \lfloor 1/f(\xi_{n}) \rfloor\}$ and let $ (x_{m})_{m < \omega}$, $ (y_{m})_{m < \omega}$ be as above. Then $(y_{m}) \bE_{f} \underline{0}$ and $\underline{0} \bE_{f} (x_{m}+y_{m})$ but $(y_{m})\not \!\!\bE_{f} (x_{m}+y_{m})$.$\bs$

\medskip

If $f$ is an arbitrary function, thus $\bE_{f}$ is not necessarily an equivalence relation, then one could consider the equivalence relation generated by $\bE_{f}$. However, it is very hard to control the properties of this generated equivalence relation by the properties of $f$, in particular we do not know how to ensure $\bE_{f}$ is Borel. Therefore, from now on, we restrict our attention to such functions $f$ for which $\bE_{f}$ is an equivalence relation. 

Despite of the general setting of Definition \ref{Csop} and Proposition \ref{gen}, in the present paper we will work only with two special cases. At some point, we will set $G= H$ to be the circle group $\CG=[0,1)$ with mod 1 addition. Then  $(H_{1})$-$(H_{3})$ obviously hold, $x^{+} = x$ $(x \in \CG)$, moreover $(R_{2}a)$ and $(R_{2}b)$ are equivalent. But mainly we will work with $G = \real$ and $H = [0,1]$; then $(H_{1})$-$(H_{3})$ hold and $x^{+} = |x|$.  Our reason for working with functions $f$ defined on $[0,1]$ instead of $\real$ is that on a smaller domain it is easier to define $f$ such that it satisfies ($R_{1}$) and $(R_{2})$. Next we show that for $\mathbf{E}_{\Id^{p}}$, this change of domain makes no difference.

\begin{lemma}\label{mege} For $1 \leq p < \infty$, let $\ell^{p}$ denote the equivalence relation defined by (\ref{defG}) with $f \colon  \real \rar \real^{+}$, $f(x) = |x|^{p}$ $(x \in \real)$. Then $\ell^{p}$ and $\ell^{p}|_{[0,1]^{\omega} \times [0,1]^{\omega}}$ are Borel equivalent.
\end{lemma}
\prf  It is obvious that $\ell^{p}|_{[0,1]^{\omega} \times [0,1]^{\omega}} \leq _{B} \ell^{p}$. To see the other direction, for every $k \in \egesz$ let $\rho_{k} \colon \real \rar [0,1]$, $$\rho_{k}(x) = \b\{ \begin{array}{ll}  1, \textrm{ if } k < \lfloor x \rfloor; \\  \{x\}, \textrm{ if } k = \lfloor x \rfloor; \\ 0, \textrm{ if } k > \lfloor x \rfloor; \end{array}\j.$$ and set $\vt \colon \real^{\omega} \rar [0,1]^{\egesz  \times \omega }$, $$\vt((x_{n})_{n < \omega}) = (\rho_{k}(x_{n}))_{k \in \egesz,n < \omega}.$$ For every $x,y \in \real$ with $| y-x |\leq 1$, we have $\rho_{k}(x) \neq \rho_{k}(y)$ only if $k =\lfloor x \rfloor$ or $k=\lfloor y \rfloor$; moreover $$|y-x| = \sum_{k \in \egesz} |\rho_{k}(y) - \rho_{k}(x) | ~(x,y \in \real).$$ Thus $$\sum_{k \in \egesz} |\rho_{k}(y) - \rho_{k}(x)|^{p} \leq |y-x|^{p}~(x,y \in \real);$$ and for $x,y \in \real$ with $| y-x |\leq 1$, $$ |y-x|^{p} \leq 2^{p}\sum_{k \in \egesz} |\rho_{k}(y) - \rho_{k}(x)|^{p}.$$ Since $(x_{n}) \ell_{p} (y_{n})$ implies $\lim_{n < \omega} |y_{n} -x_{n}|=0$, after reindexing the coordinates of its range, $\vt$ reduces $\ell^{p}$ to $\ell^{p}|_{[0,1]^{\omega} \times [0,1]^{\omega}}$, as required.$\bs$

\medskip

As we have seen already in the introduction, $\mathbf{E}_{f}$ may be an equivalence relation for a discontinuous $f$, e.g., 
for the function $f(0) = 0$, $f(x) = 1$ $(0 < x \leq 1)$ we have $\mathbf{E}_{f}$ is the equivalence relation of eventual equality on $[0,1]^{\omega}$. Following the literature, we denote this equivalence relation by $E_{1}$. In the remaining part of this section we show that $f$ is continuous in zero if and only if $E_{1} \not \leq _{B} \bE_{f}$.

\begin{theorem}\label{folyt} Let $ f \colon [0,1] \rar \real^{+}$ be a bounded Borel function such that $\bE_{f}$ is an equivalence relation. Then $f$ is continuous in zero if and only if $E_{1} \not \leq _{B} \bE_{f}$.
\end{theorem}

Before proving Theorem \ref{folyt} we show that up to Borel reducibility, requiring continuity in zero or continuity on the whole $[0,1]$ is the same condition for $\bE_{f}$.

\begin{clm}\label{f_e} Let $ f \colon [0,1] \rar \real^{+}$ be a bounded function such that $\bE_{f}$ is an equivalence relation. If $f$ is continuous in zero then there exists a continuous function $\tilde{f} \colon [0,1] \rar \real^{+}$ such that $\bE_{f} = \bE_{\tilde{f}}$.
\end{clm}

As a corollary of Theorem \ref{folyt} and Proposition \ref{f_e}, we obtain the following surprising result.

\begin{corollary}\label{fkov}
Let $ f,g \colon [0,1] \rar \real^{+}$ be bounded Borel functions such that $\bE_{f}$ and $\bE_{g}$  are equivalence relations. If $g$ is continuous and $\bE_{f} \leq _{B} \bE_{g}$ then $f$ is continuous in zero hence there is a continuous function $\tilde{f} \colon [0,1] \rar \real^{+}$ such that $\bE_{f}=\bE_{\tilde{f}}$.
\end{corollary}

We start with the proof of Proposition \ref{f_e}.

\medskip

\textbf{Proof of Proposition \ref{f_e}.} Let $C \geq 1$ be the constant of $(R_{2})$.  First we show that there exists an increasing function $\ve \colon  [0,1] \rar [0,1]$ such that $\ve(a) > 0$ for $a > 0$ and for every $x,y \in [0,1]$, \Keq\label{34} |y-x| \leq \ve(f(x)) \Rightarrow \frac{f(x)}{2C} \leq f(y) \leq 2Cf(x).\Zeq Set $$\ve(a)= \frac{1}{2}\sup\b\{y \in [0,1] \colon f(d) \leq \frac{a}{2C} \textrm{ for } 0 \leq d \leq y\j\};$$ then $\ve$ is increasing and since $f(0) = 0$ and $f$ is continuous in zero,  $\ve(a) >0$ for $a > 0$. We show (\ref{34}). By $(R_{2}a)$, $$f(y) \leq C(f(x)+f(y-x)) \leq 2Cf(x)~(0 \leq y-x \leq \ve(f(x)))$$ and $$\frac{f(x)}{2C} \leq \frac{f(x)}{C} - f(x-y) \leq f(y)~(0 \leq x-y \leq \ve(f(x)));$$ and by $(R_{2}b)$, $$\frac{f(x)}{2C} \leq \frac{f(x)}{C} - f(y-x) \leq f(y)~(0 \leq y-x \leq \ve(f(x)))$$ and $$f(y) \leq C(f(x)+f(x-y)) \leq 2Cf(x)~(0 \leq x-y \leq \ve(f(x))),$$ as required.

As a corollary of (\ref{34}), we get $U=\{x \in [0,1] \colon f(x) >0\}$ is an open set. Moreover, for every $a > 0$ the $\ve(a)$-neighborhood of $\{x \in [0,1] \colon f(x) > a\}$ is contained in $U$, i.e.\ $f$ is continuous at every point of $[0,1] \sm U = \{x \in [0,1] \colon f(x) = 0\}$. For every $x \in U$, set $$I_{x} = (x-\ve(f(x)), x+\ve(f(x))) \cap [0,1].$$ Then $I_{x} \ss U$ and $\{I_{x} \colon x \in U\}$ is an open cover of $U$. Since the covering dimension of $U$ is one, there is an open refinement $J_{x} \ss I_{x}$ $(x \in U)$ such that $\{J_{x} \colon x \in U\}$  is an open cover of $U$ of order at most two, i.e.\ for every $x \in U$, $|\{y \in U \colon x \in J_{y}\}| \leq 2$. So the function $\vp \colon U \rar 2^{\real}$, \Keq \notag \vp(x) = \bigcup_{\footnotesize \begin{array}{c}y \in U \\ x \in J_{y}\end{array}} \b[ \frac{f(y)}{2C},2Cf(y)\j] \Zeq is closed convex valued and lower semicontinuous, hence Michael's Selection Theorem \cite[Theorem 3.2 p.\ 364]{M} can be applied to have a continuous function $\tilde f \colon U \rar \real$ satisfying $\tilde f (x) \in \vp(x)$ $(x \in U)$. Since $f$ is continuous at every point of $[0,1] \sm U$, $\tilde f$ extends continuously to $[0,1]$ with $\tilde f (x) = 0$ for $x \in  [0,1] \sm U$.

For fixed $x \in [0,1]$, $x \in J_{y}$ implies $x \in I_{y}$. So by (\ref{34}), $$f(x) \in \b[ \frac{f(y)}{2C},2Cf(y)\j]\textrm{ hence } f(y) \in \b[ \frac{f(x)}{2C},2Cf(x)\j].$$ Thus $$\bigcup_{\footnotesize \begin{array}{c}y \in U \\ x \in J_{y}\end{array}} \b[ \frac{f(y)}{2C},2Cf(y)\j] \ss \b[ \frac{f(x)}{4C^2},  4C^2f(x)  \j] $$ and so $$\frac{f(x)}{4C^2} \leq \tilde f(x) \leq 4C^2f(x) ~(x \in [0,1]).$$ Therefore $\bE_{f}=\bE_{\tilde f}$, as required.$\bs$

\medskip

We close this section with the proof of Theorem \ref{folyt}. We obtain the nonreducibility of $E_{1}$ to $\bE_{f}$ for a continuous $f$ via \cite[Theorem 4.1 p.\ 238]{KL}, which says that $E_{1}$ is not reducible to any equivalence relation induced by a Polish group action. To this end, first we show that for continuous $f$, $\bE_{f}$ is essentially induced by a Polish group action. Recall that $\CG$ denotes the circle group $[0,1)$ with mod 1 addition.

\begin{lemma} \label{cont} Let $f \colon [0,1] \rar \real^{+}$ be a continuous function such that $\bE_{f}$ is an equivalence relation. Then either $f$ is identically zero or there is a continuous even function $\tilde f \colon \CG \rar \real^{+}$ such that $\tilde f(x) > 0$ for $x \neq 0$, $\bE_{\tilde f}$ is an equivalence relation and $\bE_{f} \leq _{B} \bE_{\tilde f}$.
\end{lemma}
\prf Suppose $f$ is not identically zero. We distinguish two cases. Suppose first $f(x) > 0$ for $x > 0$. Then set $$\tilde f(x) = \b\{ \begin{array}{ll} f(2x), & \textrm{ if } 0 \leq x < 1/2;  \\ f(2-2x), & \textrm{ if } 1/2 \leq x < 1.  \end{array} \j.$$ It is obvious that $\tilde f$ is continuous, even and $\tilde f(x) > 0$ for $x \neq 0$. We show that $\bE_{\tilde f}$ is an equivalence relation by verifying the conditions of Proposition \ref{ekv}. We have $(R_{1})$; since $(R_{2}a)$ implies $(R_{2}b)$, we prove only $(R_{2}a)$. Let $C$ be the constant of $(R_{2})$ for $f$. If $x \in [1/4,3/4]$ or $y \in [1/4,3/4]$ then $$\tilde f (x+y) \leq \frac{\max \tilde f}{ \min \tilde f|_{[1/4,3/4]}} (\tilde f(x) + \tilde f(y)).$$ If $x,y \in[0,1/4]$ or $x,y \in[3/4,1)$ then by $(R_{2}a)$ for $f$,  $(R_{2}a)$ holds for $\tilde f$ with $C$. Finally if exactly one of $x$ and $y$ is in $[0,1/4]$ and $[3/4,1)$ then by $(R_{2}b)$ for $f$,  $(R_{2}a)$ holds for $\tilde f$ with $C$.

 Also, $\vt \colon [0,1]^{\omega} \rar \CG^{\omega}$, $\vt((x_{n})_{n < \omega}) = (x_{n}/2)_{n < \omega}$ is a reduction of $\bE_{f}$ to $\bE_{\tilde f}$, so the proof of first case is complete.

In the second case, suppose $f(x) = 0$ for some $x \in (0,1]$. By $(R_{2})$, the nonempty set $\{x \in [0,1] \colon f(x) = 0\}$ is closed under 
additions that are in $[0,1]$. Hence by the continuity of $f$, $x^{\star} = \inf\{x \in (0,1] \colon f(x) = 0\}$ satisfies $x^{\star} > 0$ and $f(x^{\star}) = 0$.

Set $$\tilde f(x) = \b\{ \begin{array}{ll} f(xx^{\star}), & \textrm{ if } 0 \leq x < 1/2;  \\ f((1-x)x^{\star}), & \textrm{ if } 1/2 \leq x < 1.  \end{array} \j.$$ It is obvious that $\tilde f$ is continuous, even and $\tilde f(x) > 0$ for $x \neq 0$. Similarly to the previous case, we get $\bE_{\tilde f}$ is an equivalence relation by distinguishing several cases.  If $x \in [1/4,3/4]$ or $y \in [1/4,3/4]$ then $$\tilde f (x+y) \leq \frac{\max \tilde f}{ \min \tilde f|_{[1/4,3/4]}} (\tilde f(x) + \tilde f(y)).$$ If $x,y \in[0,1/4]$  then by $(R_{2}a)$ for $f$,  $(R_{2}a)$ holds for $\tilde f$ with $C$. If $x,y \in[3/4,1)$ then again by  $(R_{2}a)$ for $f$, $$\tilde f (x+y) = f(2  x^{\star} - (x+y) x^{\star}) \lesssim  f( x^{\star} - x x^{\star}) + f( x^{\star} - y x^{\star}) = \tilde f (x) + \tilde f (y). $$ Finally if exactly one of $x$ and $y$ is in $[0,1/4]$ and $[3/4,1)$ then by $(R_{2}b)$ for $f$,  $(R_{2}a)$ holds for $\tilde f$ with $C$. 

For every $x \in [0,1]$, let $\langle x\rangle = x/x^{\star}-\lfloor x/x^{\star} \rfloor $. We show that 
$\vt \colon [0,1]^{\omega} \rar \CG^{\omega}$, $\vt((x_{n})_{n < \omega}) = (\langle x_{n} \rangle)_{n < \omega}$ is a reduction of $\bE_{f}$ to $\bE_{\tilde f}$. For every $0\leq x \leq y \leq 1$, with $k =  \b\lfloor y/x^{\star} \j\rfloor - \b\lfloor x/x^{\star} \j\rfloor$ we have $\langle y \rangle -\langle x \rangle = y/x^{\star} - x/x^{\star} - k$, so   $$ \tilde f \b(\langle y \rangle -\langle x \rangle\j) = \b\{ \begin{array}{ll} f\b(y-x -k x^{\star} \j), & \textrm{ if } 0 \leq y-x -k x^{\star} < x^{\star}/2;  \\ f(-y+x + k x^{\star}), & \textrm{ if } 0 \leq -y+x+k x^{\star} < x^{\star}/2; \\ f\b(x^{\star} -y + x +k x^{\star} \j), & \textrm{ if } x^{\star}/2 \leq y-x -k x^{\star} < x^{\star};  \\ f(x^{\star} +y-x-k x^{\star}), & \textrm{ if } x^{\star}/2 \leq -y+x+k x^{\star} < x^{\star}.  \end{array} \j.$$ For $l = k$ or $l = k\pm 1$, in any of the cases where applicable, by $(R_{2})$ we have $$f(y-x) \lesssim f(y-x - lx^{\star}) + f(lx^{\star}), ~ f(y-x)\lesssim f(lx^{\star}) + f(lx^{\star} -y+x),$$ $$f(y-x - lx^{\star})\lesssim f(y-x) + f(lx^{\star}), ~f(lx^{\star} -y+x) \lesssim f(lx^{\star}) + f(y-x).$$ So $f(y-x) \approx \tilde f (\langle y \rangle - \langle x\rangle )$ follows from $f(lx^{\star}) = 0$. This implies that $\vt$ is a reduction, so  the proof is complete.$\bs$

\medskip

In the next lemma, for an $\tilde f$ as in Lemma \ref{cont}, we find a Polish group action inducing $\bE_{\tilde f}$.

\begin{definition}\rm\label{No} Let $f \colon H \rar \real^{+}$ be an arbitrary function. For every $x  = (x_{n})_{n < \omega} \in H^{\omega}$ and $I \ss \omega$ we set $$\|x \|_{f} = \sum_{n < \omega}f(x_{n}), ~\|x|_{I} \|_{f} = \sum_{n \in I}f(x_{n}).$$ We define $\mc{N}_{f}= \{x\in H^{\omega} \colon \|x\|_{f} < \infty\}$.
\end{definition}

\begin{lemma}\label{idea} Let $f \colon \CG \rar \real^{+}$ be a continuous even function such that $f(x) > 0$ for $x \neq 0$ and $\bE_{f}$ is an equivalence relation. 
\ben
\item\label{idea1} There is a unique topology $\tau_{f}$ on $\mc{N}_{f}$ such that for every $x \in\mc{N}_{f}$, the sets $$B(x, \ve) = \{y \in \mc{N}_{f} \colon \|y-x\|_{f} < \ve\}~(\ve > 0)$$ form a neighborhood base at $x$. This topology is regular, second countable and refines the topology inherited from $\CG^{\omega}$.

\item\label{idea2} With $\tau_{f}$, $(\mc{N}_{f},+)$ is a Polish group. The natural action of $\mc{N}_{f}$ on $\CG^{\omega}$ is continuous, and the equivalence relation induced by this action is $\bE_{f}$.
\een
\end{lemma}
\prf For \ref{idea1}, we show that for every $x \in \mc{N}_{f}$, $\ve > 0$ and $y \in B(x, \ve)$ there is a $\delta > 0$ such that  $B(y,\delta) \ss B(x,\ve)$; once this done, the first part of the statement follows from elementary topology (see e.g.\ \cite{C}). Let $C \geq 1$ be the constant of $(R_{2})$, fix $x \in \mc{N}_{f}$, $\ve > 0$ and $y \in B(x, \ve)$. Let $n < \omega$ be such that $$\|(y-x)|_{\omega\sm n} \|_{f} < \frac{\ve - \|y-x\|_{f}}{3C}.$$ Let $\delta > 0$ satisfy $\delta < (\ve-\|y-x\|_{f})/(3C)$, and such that for every $i < n$ and $z_{i} \in [0,1]$, $f(z_{i} - y_{i}) < \delta$ implies $$|f(z_{i}-x_{i}) - f(y_{i}-x_{i})| < \frac{\ve-\|y-x\|_{f}}{3n};$$ such a $\delta$ exists by the continuity of $f$ and by $f(x) > 0$ for  $x \neq 0$.

Let $z \in B(y,\delta)$; then by $(R_{2})$, \begin{multline}\notag \|z-x\|_{f} = \|(z-x)|_{n}\|_{f} + \|(z-x)|_{\omega \sm n}\|_{f} <  \\ \|(y-x)|_{n}\|_{f} + n \frac{\ve-\|y-x\|_{f}}{3n} + C(\|(z-y)|_{\omega \sm n}\|_{f} + \|(y-x)|_{\omega \sm n}\|_{f}) < \\ \|y-x\|_{f} + \frac{\ve-\|y-x\|_{f}}{3} + \frac{\ve-\|y-x\|_{f}}{3} + \frac{\ve-\|y-x\|_{f}}{3} = \ve,
\end{multline} as required.

Since $f(x) > 0$ for $x \neq 0$, $\tau_{f}$ refines the topology inherited from $\CG^{\omega}$. The countable set of eventually zero rational sequences shows separability and hence second countability.
To see regularity, let $F \ss (\mc{N}_{f}, \tau_{f})$ be a closed set and take $x \notin F$. Then $\gamma = \inf\{\|y-x\|_{f} \colon y \in F\} > 0$. By $(R_{2})$, $B(x,\gamma/(2C)) \cap B(F,\gamma/(2C)) = \es$, as required.

For \ref{idea2}, first we show $(\mc{N}_{f},+)$ is a topological group. Let $x,y \in \mc{N}_{f}$ and $\gamma > 0$. By $(R_{2})$, $B(x,\gamma/2C)+B(y,\gamma/2C) \ss B(x+y,\gamma)$, so addition is continuous. The continuity of the inverse operation is obvious, so the statement follows.

Next we show $(\mc{N}_{f}, \tau_{f})$ is strong Choquet (for the definition and notation see \cite[Section 8.D p.\ 44]{K}). The closed balls $\overline B(x,\ve) = \{y \in \mc{N}_{f} \colon \|y-x\|_{f}  \leq \ve\}$ are closed in $\CG^{\omega}$, thus every $\| \cdot\|_{f}$-Cauchy sequence is convergent in $\mc{N}_{f}$. If player $I$ plays $(x_{n}, U_{n})_{n < \omega}$, a winning strategy for player $II$ is to choose $V_{n}=B(x_{n}, \gamma_{n})$ such that $\overline B(x_{n}, \gamma_{n}) \ss U_{n}$ and $\gamma_{n} \leq 1/2^{n}$ $(n < \omega)$. So  $(\mc{N}_{f}, \tau_{f})$ is  strong Choquet, hence Polish by Choquet's Theorem (see e.g.\ \cite[(8.18) Theorem p.\ 45]{K}). 

The continuity of the action of $\mc{N}_{f}$ on $\CG^{\omega}$ follows from the fact that $\tau_{f}$ refines the topology inherited from $\CG^{\omega}$. It is obvious that the equivalence relation induced by this action is $\bE_{f}$, so the proof complete.$\bs$

\medskip

\begin{definition}\rm\label{V} For a topological space $X$ and $G \ss X$, we set $$V(G)= \bigcup\{ U \ss X \colon U \textrm{ is open, } G \cap U \textrm{ is comeager in } U\}.$$
\end{definition}

\medskip
 
 The next lemma is a folklore result on the existence of a perfect set with special distance set.

\begin{lemma} \label{dis} Let $G \ss [0,1]$ be a Borel set such that zero is adherent to $V(G)$. Then there exists a nonempty perfect set $P \ss [0,1]$ such that \Keq\label{DI}\{|y-x| \colon x,y \in P\} \ss G \cup \{0\}.\Zeq
\end{lemma}
\prf By passing to a subset, we can assume that $G$ is a comeager $G_{\delta}$ subset of $V(G)$. Set $\tilde G = G \cup (-G)\cup \{0\}$ and let $d_{\tilde G}$ be the metric on $\tilde G$ for which $(\tilde G,d_{\tilde G})$ is a Polish space with the topology inherited from $[-1,1]$ (see e.g.\ \cite[(3.11) Theorem p.\ 17]{K}). We construct inductively a sequence $(x_{n})_{n < \omega} \ss [0,1]$ with the following properties:
\ben
\item\label{en0} for every $n < \omega$, $x_{n+1} < x_{n}/2$;
\item\label{en1} for every $ s \in \{-1,0,+1\}^{<\omega}$, $\sum_{i <|s|} s(i)x_{i} \in \tilde G$; 
\item\label{en2} for every $ s \in \{-1,0,+1\}^{<\omega}\sm \{\es\}$, $d_{\tilde G}(\sum_{i <|s|-1} s(i)x_{i}, \sum_{i <|s|} s(i)x_{i}) \leq 1/2^{|s|}$.
\een
Let $x_{0} \in G$ be arbitrary. Let $0 < n < \omega$ and suppose $x_{i}$ $(i < n)$ are defined such that \ref{en1} and \ref{en2} hold for every $s \in  \{-1,0,+1\}^{\leq n}$. By  \ref{en1}, if $s \in  \{-1,0,+1\}^{n}$ and $\sum_{i <n} s(i)x_{i}  \neq 0$ then $\tilde G$ is comeager in a neighborhood of $\sum_{i <n} s(i)x_{i}$. Since zero is adherent to $V(G)$, by the Baire Category Theorem we can pick $x_{n} \in G$ sufficiently close to zero such that \ref{en0} holds; and for every $s \in  \{-1,0,+1\}^{n+1}$ with $\sum_{i <n} s(i)x_{i}  \neq 0$ we have $\sum_{i <n+1} s(i)x_{i}  \in \tilde G$, hence by $x_{n} \in G$, $\sum_{i <n+1} s(i)x_{i}  \in \tilde G$ for every $s \in  \{-1,0,+1\}^{n+1}$; and in addition  \ref{en2} holds. This completes the inductive step.

We show $$P = \b\{\sum_{n < \omega}\sigma(n)x_{n} \colon \sigma \in 2^{\omega}\j\}$$ fulfills the requirements. By \ref{en2}, for every $\sigma \in  \{-1,0,+1\}^{\omega}$, $(\sum_{i <n} \sigma(i)x_{i})_{n < \omega}$ is a Cauchy sequence in $\tilde G$, so $\sum_{n <\omega} \sigma(n)x_{n} \in \tilde G$.  In particular $P \ss G\cup \{0\}$.

Let $x,x' \in P$, $x=\sum_{n < \omega}\sigma(n)x_{n}$ and $x' = \sum_{n < \omega}\sigma'(n)x_{n}$ with $\sigma, \sigma' \in 2^{\omega}$, $\sigma \neq \sigma'$; say for the first $n < \omega$ with $\sigma(n) \neq \sigma'(n)$ we have $\sigma(n) = 0$, $\sigma'(n) = 1$. Then  for $\delta \in \{-1,0,+1\}^{\omega}$, $\delta(n) = \sigma'(n) - \sigma(n)$ $(n < \omega)$ we have $$|x' - x| = x' - x = \sum_{ n < \omega} \delta (n) x_{n} \in \tilde G, $$ moreover by \ref{en0}, $x' - x > 0$ i.e.\ $x' - x \in G$. Thus $P$ is a nonempty perfect set and satisfies (\ref{DI}), which completes the proof.$\bs$

\medskip

The last lemma points out a property of an $f$  discontinuous in zero. 

\begin{lemma} \label{dis1} Let $ f \colon [0,1] \rar \real^{+}$ be a bounded Borel function such that $\bE_{f}$ is an equivalence relation. If $f$ is not continuous in zero then there exists an $a > 0$ such that  $G = \{x \in [0,1] \colon f(x) > a\}$ satisfies the condition of Lemma \ref{dis}, i.e.\ zero is adherent to $V(G)$.
\end{lemma}
\prf Let $C$ be the constant of $(R_{2})$. Since $f$ in not continuous in zero, there exists an $a > 0$ such that zero is adherent to $\{x \in [0,1] \colon f(x) > 2Ca\}$. If for every $x \in [0,1]$ with $f(x) > 2Ca$, $x$ is adherent to $$\bigcup \{U \ss (x,1) \colon U \textrm{ is open, } \{y \in U \colon f(y) > a\} \textrm{ is comeager in } U\}$$ then the statement follows. If not, by $f$ being Borel, there is an $x \in [0,1]$ with $f(x) > 2Ca$  and a $\delta > 0$ such that $$Y = \{y \in (x, x+\delta) \colon f(y) \leq  a\}$$ is comeager in $(x, x+\delta)$. Since $f(x) > 2Ca$, by $(R_{2}b)$ we have $f(y-x) >a$ whenever  $y \in Y$. Hence $\{x \in [0,1] \colon f(x) > a\}$ is comeager in $(0,\delta)$, which finishes the proof.$\bs$

\medskip

\textbf{Proof of Theorem \ref{folyt}.} Suppose first $f$ is not continuous in zero. By Lemma \ref{dis} and   Lemma \ref{dis1}, there is an $a > 0$ and a nonempty perfect set $P \ss [0,1]$ such that $f(|y - x|) > a$ for every $x,y \in P$, $x \neq y$. Thus $\bE_{f}$ restricted to $P^{\omega}$ is $E_{1}$.

Suppose now that $f$ is continuous in zero. By Proposition \ref{f_e}, we can assume $f$ is continuous on $[0,1]$. If $f \equiv 0$ then $E_{1} \not \leq _{B} \bE_{f}$ is obvious. Else by Lemma \ref{cont} and Lemma \ref{idea}, $\bE_{f} \leq _{B}\bE_{\tilde f}$ where $\bE_{\tilde f}$ is induced by a Polish group action. Hence $E_{1} \not \leq _{B} \bE_{\tilde f}$ by \cite[Thorem 4.1 p.\ 238]{KL} and \cite{Ke}; in particular $E_{1} \not \leq _{B} \bE_{f}$. This completes the proof.$\bs$

\section{Reducibility results} \label{red}

In the remaining part of the paper, in most cases, we restrict our attention to equivalence relations $\bE_{f}$ where $ f \colon [0,1] \rar \real^{+}$ is a \emph{continuous} function. As we have seen in Proposition \ref{f_e}, requiring continuity on $[0,1]$ and continuity in zero for $f$ are equivalent, and by Theorem \ref{folyt}, for Borel $f$  it is the necessary and sufficient condition to have $E_{1} \not \leq_{B}  \bE_{f}$. This assumption is acceptable to us since we aim to study equivalence relations $\bE_{f}$ for which $\bE_{f} \leq _{B} \bE_{\Id^{q}}$ for some $1 \leq q < \infty$.

The main  restriction, in addition to $(R_{1})+(R_{2})$, we impose  in the sequel on the function $f$ is formulated in the following definition.

\begin{definition}\label{essi}\rm Let $(R ,\leq)$ be an ordered set and $f \colon R \rar \real^{+}$ be a function. We say $f$ is \emph{essentially increasing} if for some $C\geq 1$,  $\forall x, y \in R$ $(x \leq y\Rightarrow f(x) \leq Cf(y))$. Similarly, $f$ is \emph{essentially decreasing} if for some $C\geq 1$, $\forall x, y \in R$ $(x \leq y \Rightarrow C f(x) \geq f(y))$. 
\end{definition}

\begin{lemma}\label{EINC} With the notation of Definition \ref{essi}, $f$ is essentially increasing (resp.\ essentially decreasing) if and only if there is an increasing (resp.\ decreasing) function  $\tilde f$ such that $\tilde f \approx f$. 
\end{lemma}
\prf If $f$ is essentially increasing, set $\tilde f \colon R \rar \real^{+}$, $\tilde f (x) = \sup \{f(y) \colon y \leq x\}$. Then $\tilde f$ is increasing and $f(x) \leq \tilde f(x) \leq C f(x)$ $(x \in R)$.  If  $f$ is essentially decreasing, let  $ \tilde f (x) = \inf \{f(y) \colon y \leq x\}$; then,  as above, $\tilde f \approx  f$. The other directions are obvious, so the proof is complete.$\bs$

\medskip
 
 We remark that for $R = [0,1]$ or $R = (0,1]$, by its definition above, $\tilde f$ is continuous if $f$ is so. By $\tilde f \approx f$, $\tilde f$ and $f$ have the same asymptotic behavior in 0.
 
In this section we prove the following two theorems.

\begin{theorem}\label{ala} Let $1 \leq \alpha < \infty$ and let $\psi \colon (0,1] \rar (0,+\infty)$ be an essentially decreasing continuous function such that $\Id^{\alpha} \psi$ is bounded and for every $\delta > 0$, $\liminf_{x \rar +0} x^{\delta} \psi(x) = 0$. Set $g(x) = x^{\alpha}\psi(x)$ for $0 < x \leq 1$ and $g(0)=0$. Suppose $\bE_{g}$ is an equivalence relation. Then $\bE_{\mathrm{Id}^{\alpha}} \leq _{B} \bE_{g}$.
\end{theorem}

\begin{theorem}\label{fole} Let $f,g \colon [0,1] \rar \real^{+}$ be continuous essentially increasing functions such that $\bE_{f}$ and $\bE_{g}$ are equivalence relations. Suppose there exists a function $\kappa \colon \{1/2^i \colon i < \omega\} \rar [0,1]$ satisfying the recursion \Keq\label{U2}f(1)=g(\kappa(1)),~f(1/2^{n})=\sum_{i=0}^{n} g(\kappa(1/2^{i})/2^{n-i})~(0 < n < \omega)\Zeq such that for some $L \geq 1$,  \Keq\label{U3}\sum_{i=n}^{\infty} g(\kappa(1/2^{i})) \leq L \sum_{i=0}^{n}  g(\kappa(1/2^{i})/2^{n-i}) ~(n < \omega)\Zeq and \Keq\label{star} \kappa(1/2^{n}) \leq L \cdot \max \{\kappa(1/2^{i})/2^{n-i}\colon i < n\}~(n < \omega).\Zeq Then  $\bE_{f} \leq _{B} \bE_{g}$.
\end{theorem}

Theorem \ref{ala} illustrates, e.g.\ by choosing $\psi(x) = 1-\log(x)$ $(0 < x \leq 1)$, that reducibility among the $\bE_{f}$s is \emph{not} characterized by the growth order of the $f$s. Theorem \ref{fole} is a stronger version of \cite[Theorem 1.1 p.\ 1836]{DH}, but we admit that our improvement is of technical nature. However, in Section \ref{con} it will allow us to show the reducibility among $\bE_{f}$s for new families of $f$s.

These results neither give a complete description of the reducibility between the equivalence relations $\bE_{f}$ nor are optimal. Nevertheless, we note that in
Theorem \ref{ala}, $\mathrm{Id}^{\alpha}$ cannot be replaced by an arbitrary ``nice" function: as we will see, e.g.\ $\bE_{\mathrm{Id}^{\alpha}} < _{B} \bE_{\mathrm{Id}^{\alpha}/(1-\log)}$. Also, the condition $\psi$ is decreasing cannot be left out: e.g.\ we need the techniques of Theorem \ref{fole} in order to treat the $\psi(x) = x$ case, i.e.\ to show $\bE_{\mathrm{Id}^{\alpha}} \leq _{B} \bE_{\mathrm{Id}^{\alpha+1}}$. We comment on the optimality of Theorem \ref{fole} after its proof.

We start with a technical lemma.

\begin{lemma} \label{disz} Let $f,g \colon [0,1] \rar \real^{+}$ be continuous functions such that $\bE_{f}$, $\bE_{g}$ are equivalence relations. Suppose there exists $K > 0$ and $I \in [\omega]^{\omega}$ such that for every $n \in I$ there is a mapping $\vt_{n} \colon \{i/n \colon 0 \leq i \leq n\} \rar [0,1]^{\omega}$ satisfying \Keq\label{U1}\frac{1}{K} f((j-i)/n) \leq \|\vt_{n}(j/n)-\vt_{n}(i/n)\|_{g} \leq K f((j-i)/n)~(0 \leq i < j \leq n).\Zeq Then $\bE_{f} \leq _{B} \bE_{g}$.
\end{lemma}
\prf For $x \in [0,1]$ and $0 < n  < \omega$ set $[x]_{n} = \max\{i/n\colon i/n \leq x, ~0 \leq i \leq n\}$. Since $f$ is uniformly continuous on $[0,1]$, for every $k < \omega$ there is an $n_{k} \in I$ such that $|f(x) - f([x]_{n_{k}})| \leq 1/2^{k}$ $(x \in [0,1])$. We show that $\vt \colon [0,1]^{\omega} \rar [0,1]^{\omega \cdot \omega}$, $$\vt((x_{k})_{k < \omega}) = (\vt_{n_{k}}([x_{k}]_{n_{k}}))_{k < \omega},$$ after reindexing the coordinates of the range, is a Borel reduction of $\bE_{f}$ to $\bE_{g}$.

Let $(x_{k})_{k < \omega}, (y_{k})_{k < \omega} \in [0,1]^{\omega}$. We have $$[|y_{k} - x_{k}|]_{n_{k}} \leq |[y_{k}]_{n_{k}} - [x_{k}]_{n_{k}}| \leq [|y_{k} - x_{k}|]_{n_{k}}  + 1/n_{k}~(k < \omega).$$ So by the choice of $n_{k}$, $|f(|y_{k} - x_{k}|) -f([|y_{k} - x_{k}|]_{n_{k}} )| \leq 1/2^{k}$  and  $|f([|y_{k} - x_{k}|]_{n_{k}} )- f(|[y_{k}]_{n_{k}} - [x_{k}]_{n_{k}}|)| \leq 1/2^{k}$, thus  $$|f(|y_{k} - x_{k}|) - f(|[y_{k}]_{n_{k}} - [x_{k}]_{n_{k}}| ) |\leq 2/2^{k}~(k < \omega).$$ By (\ref{U1}), $\|\vt_{n_{k}}([y_{k}]_{n_{k}}) - \vt_{n_{k}}([x_{k}]_{n_{k}})\|_{g} \approx f(|[y_{k}]_{n_{k}} - [x_{k}]_{n_{k}}| )$ $(k < \omega)$, so the statement follows.$\bs$

\medskip

\textbf{Proof of Theorem \ref{ala}.} For some $B \geq 1$, let $x^{\alpha}\psi(x) \leq B$ $(0 <x \leq 1)$. We find a $K > 0$ such that for every $0<n < \omega$ there exist $M < \omega$ and $0 < \mu \leq 1$ such that for every $0 \leq i < j \leq n$, \Keq\label{A1} \frac{1}{K}\b(\frac{j-i}{n}\j)^{\alpha} \leq M \b(\frac{j-i}{n}\mu\j)^{\alpha} \psi\b(\frac{j-i}{n}\mu\j)\leq K\b(\frac{j-i}{n}\j)^{\alpha}.\Zeq Once this done, the conditions of Lemma \ref{disz} are satisfied by the mapping $\vt_{n} \colon \{i/n \colon 0 \leq i \leq n\} \rar [0,1]^{\omega}$, $$\vt_{n}(i/n) = (\underbrace{i\mu/n, \dots, i\mu/n}_{M},0, \dots).$$ Observe that (\ref{A1}) is equivalent to $$1/K\leq M \mu^{\alpha} \psi((j-i)\mu/n) \leq K~(0 \leq i < j \leq n).$$ Since $\psi$ is essentially decreasing, it is enough to have $1/2 \leq M\mu^{\alpha}\psi(\mu)$ and $M\mu^{\alpha}\psi(\mu/n) \leq 2B$. We will find a $0 < \mu \leq 1$ satisfying $ \psi(\mu/n) \leq 2\psi(\mu)$. Then by choosing $M$ to be minimal such that $1/2 \leq M\mu^{\alpha}\psi(\mu)$, by $\mu^{\alpha}\psi(\mu) \leq B$ and $B \geq 1$ we have $M\mu^{\alpha}\psi(\mu/n) \leq 2 M\mu^{\alpha}\psi(\mu)  \leq 2B$, so  we fulfilled the requirements.

Suppose such a $\mu$ does not exist, i.e.\ $\psi(\mu/n) > 2 \psi(\mu)$ $(0 < \mu \leq 1)$. Then for every $k < \omega$ and  $\mu \in [1/n,1]$, $\psi(n^{-k}\mu) \geq 2^k \psi(\mu)$. We have $x = n^{-k}\mu$ runs over $(0,1]$ as $(k,\mu)$ runs over $\omega \times [1/n,1]$. So since $\psi$ is essentially decreasing, with $\delta = \log(2)/\log(n)$ we have $\psi(x)x^{\delta} \gtrsim 1/n^{\delta}\psi(1) >0$ $(0 < x \leq 1)$. This contradicts $\liminf_{x \rar +0} x^{\delta} \psi(x) = 0$, so the proof is complete.$\bs$

\medskip

\textbf{Proof of Theorem \ref{fole}.} Let $n < \omega$ be fixed. For $0 < l \leq 2^{n}$ let $r(l) \leq n$, $s(l) < \omega$ be such that $l/2^{n} = s(l)/2^{r(l)}$ and $s(l)$ is odd. With $\PR_{l}x$ standing for the $l^{\textrm{th}}$ coordinate of $x \in [0,1]^{2^{n}}$, for every $0 \leq i \leq 2^{n}$ we define $\vt(i/2^{n})$ by  \Kem\label{DD1}\PR_{l}\vt(i/2^{n}) = \DS \b(1-2^{r(l)}\b|\frac{i}{2^{n}}-\frac{l}{2^{n}}\j|\j)\kappa(1/2^{r(l)}) \textrm{ if } \\ l > 0 \textrm{ and } \b|\frac{i}{2^{n}}-\frac{l}{2^{n}}\j| \leq 1/2^{r(l)},\end{multline} else let $\PR_{l}\vt(i/2^{n}) =0$. We show (\ref{U1}) holds for $\vt_{2^{n}} = \vt$. 

Let $0 \leq i < j \leq 2^{n}$ be arbitrary. Let $m \leq n$ be minimal such that for some $e < 2^{m}$ we have $$\frac{i}{2^{n}} \leq \frac{e}{2^{m}} < \frac{(e+1)}{2^{m}} \leq \frac{j}{2^{n}}.$$ 

We distinguish several cases. 

Suppose first $i/2^{n} =e/2^{m}$ and $j/2^{n} =(e+1)/2^{m}$. For every $k \leq m$ there is exactly one $l$ with $r(l)=k$ such that $$|e/2^{m}-l/2^{n}| \leq 1/2^{k} ~and~ |(e+1)/2^{m}-l/2^{n}| \leq 1/2^{k};$$ and for this $l$, by (\ref{DD1}), $$\b|\PR_{l}(\vt((e+1)/2^{m}) - \vt(e/2^{m}))\j| = \kappa(1/2^{k})/2^{m-k}.$$ All the other coordinates of $\vt((e+1)/2^{m}) $ and $\vt(e/2^{m})$ are zero so by (\ref{U2}), \Keq\label{ezis}\|\vt((e+1)/2^{m}) - \vt(e/2^{m}))\|_{g} = \sum_{k=0}^{m} g(\kappa(1/2^{k})/2^{m-k}) = f(1/2^{m}),\Zeq i.e.\ (\ref{U1}) holds with $K=1$.

Next suppose $e$ is even, $i/2^{n} =e/2^{m}$ and $(e+1)/2^{m}<j/2^{n}$; then we have $m \geq 1$. Observe that by the choice of $m$ we have  $j/2^{n}< (e+2)/2^{m}$. For every $k < m$ there is exactly one $l$ with $r(l)=k$ such that $$|e/2^{m}-l/2^{n}| \leq 1/2^{k} ~and~|j/2^{n}-l/2^{n}| \leq 1/2^{k};$$ and for this $l$, $l/2^{n} \notin (e/2^{m},j/2^{n})$. So by (\ref{DD1}), $$\DS \frac{\kappa(1/2^{k})}{2^{m-k}} \leq \b| \PR_{l} (\vt(j/2^{n}) - \vt(e/2^{m}))\j| \leq \DS 2\frac{\kappa(1/2^{k})}{2^{m-k}} ~(0 \leq k < m).$$ Since $e$ is even, $\vt(e/2^{m})$ has no other nonzero coordinates. For every $m \leq k \leq n$ there is exactly one $l$ with $r(l)=k$ such that $ |j/2^{n}-l/2^{n}| \leq 1/2^{k},$ and for this $l$, $\PR_{l}(\vt(j/2^{n})) \leq \kappa(1/2^{k})$. Since $g$ is essentially increasing, we have \begin{multline}\label{U6} \sum_{k=0}^{m-1} g\b(\frac{\kappa(1/2^{k})}{2^{m-k}}\j) \lesssim \|\vt(j/2^{n}) - \vt(e/2^{m}))\|_{g} \lesssim \\ \sum_{k=0}^{m-1} g\b(2\frac{\kappa(1/2^{k})}{2^{m-k}}\j)+ \sum_{k=m}^{n} g(\kappa(1/2^{k})).\end{multline}
 By $(R_{2})$,  \Keq\label{U5}g\b(2\frac{\kappa(1/2^{k})}{2^{m-k}}\j) \lesssim g\b(\frac{\kappa(1/2^{k})}{2^{m-k}}\j) ~(0 \leq k < m).\Zeq By (\ref{U2}) and since $f$ is essentially increasing, $$f\b(\frac{j}{2^n} - \frac{e}{2^m}\j) \lesssim f\b(\frac{1}{2^{m-1}}\j) = \sum_{k=0}^{m-1} g\b(\frac{\kappa(1/2^{k})}{2^{m-1-k}}\j),$$ so by (\ref{U6}), \Keq\label{U11}f\b(\frac{j}{2^n} - \frac{e}{2^m}\j)  \lesssim \|\vt(j/2^{n}) - \vt(e/2^{m}))\|_{g}.\Zeq
By (\ref{U3}) and (\ref{U5}), the right hand side of (\ref{U6}) is $\lesssim \sum_{k=0}^{m} g\b(\kappa(1/2^{k})/2^{m-k}\j)$, so since $f$ is essentially increasing, \begin{multline}\label{U10}\|\vt(j/2^{n}) - \vt(e/2^{m}))\|_{g} \lesssim  \\ \sum_{k=0}^{m} g\b(\frac{\kappa(1/2^{k})}{2^{m-k}}\j)  = f(1/2^{m}) \lesssim f\b(\frac{j}{2^n} - \frac{e}{2^m}\j) .\end{multline} The case $e+1$ is even, $i/2^{n} < e/2^{m}$ and $j/2^{n} = (e+1)/2^{m}$ can be treated by an analogous argument.

Suppose now $e$ is even, $i/2^{n} < e/2^{m}$ and $(e+1)/2^{m} \leq j/2^{n}$; then we have $m \geq 2$. By $(R_{2})$, (\ref{U10}), and also by (\ref{ezis}) if $j/2^{n}=(e+1)/2^{m} $,  \begin{multline}\notag\|\vt(j/2^{n}) - \vt(i/2^{n}))\|_{g} \lesssim \\ \|\vt(j/2^{n}) - \vt(e/2^{m}))\|_{g} +\|\vt(e/2^{m}) - \vt(i/2^{n}))\|_{g}  \lesssim \\ \b(f\b(\frac{j}{2^n} - \frac{e}{2^m}\j)+ f\b(\frac{e}{2^m} - \frac{i}{2^n}\j)\j) \lesssim f\b(\frac{j-i}{2^n} \j).\end{multline} To have a lower bound, observe that for every $k < m-1$ there is exactly one $l$ with $r(l)=k$ such that $$|i/2^{n}-l/2^{n}| \leq 1/2^{k}~and~ |j/2^{n}-l/2^{n}| \leq 1/2^{k}.$$ For this $l$, $l/2^{n} \notin (i/2^{n},j/2^{n})$ if $l/2^{n} \neq e/2^{m}$, i.e.\ if $k \neq k_{0} = r(2^{n-m}e)$. So by (\ref{DD1}), $$\DS \frac{\kappa(1/2^{k})}{2^{m-k}} \leq \b|\PR_{l} (\vt(j/2^{n}) - \vt(i/2^{n}))\j|~(k < m-1, ~k \neq k_{0})).$$ By (\ref{star}), $(R_{2})$ and since $g$ is essentially increasing, $$g\b(\frac{\kappa(1/2^{k_{0}})}{2^{m-k_{0}}}\j) \lesssim \max_{i < k_{0}}g\b(\frac{\kappa(1/2^{i})}{2^{m-i}}\j).$$ So by $(R_{2})$ and since $g$ is essentially increasing,  \begin{multline}\notag f(1/2^{m-2}) = \sum_{k=0}^{m-2} g\b(\frac{\kappa(1/2^{k})}{2^{m-2-k}}\j) \lesssim\sum_{k=0}^{m-2} g\b(\frac{\kappa(1/2^{k})}{2^{m-k}}\j)  \lesssim \\ \sum\b\{g\b(\frac{\kappa(1/2^{k})}{2^{m-k}}\j) \colon k < m-1,~k \neq k_{0} \j\} \lesssim \|\vt(j/2^{n}) - \vt(i/2^{n}))\|_{g}.\end{multline}  By the choice of $m$ we have $(j-i)/2^{n}< 4/2^{m}$. So  since $f$ is essentially increasing, $f((j-i)/2^n) \lesssim f(4/2^m) = f(1/2^{m-2})$; thus  $f\b((j-i)/2^n \j) \lesssim \|\vt(j/2^{n}) - \vt(i/2^{n})\|_{g}$. 

The case $e+1$ is even, $i/2^{n} \leq e/2^{m}$ and $(e+1)/2^{m}<j/2^{n}$ follows similarly, so the proof is complete.$\bs$

\medskip

The assumptions of Theorem \ref{fole} are not necessary, they merely make possible to imitate the construction in  the proof of \cite[Theorem 1.1 p.\ 1836]{DH}. We note however that the problem of characterizing whether $\{i/2^n \colon 0 \leq i \leq 2^n\}$ endowed with the $\|\cdot\|_{f}$-distance Lipschitz embeds into $[0,1]^{\omega}$ endowed with the $\|\cdot\|_{g}$-distance is very hard even if the distances   $\|\cdot\|_{f}$ and $\|\cdot\|_{g}$ can be related to norms (see e.g.\ \cite{Ma} and the references therein). So it is unlikely that there is a simple characterization of reducibility among $\bE_{f}$s using the approach of Lemma \ref{disz}.

\section{Nonreducibility results} \label{nonred}

In this section we improve \cite[Theorem 2.2 p.\ 1840]{DH} in order to obtain nonreducibility results for a wider class of $\bE_{f}$s, as follows. 

\begin{theorem}\label{Nr} Let $1 \leq \alpha < \infty$ and let $\vp, \psi \colon [0,1] \rar [0,+\infty)$ be continuous functions. Set $f=\Id^{\alpha}\vp$, $g(x) = \Id^{\alpha}\psi$ and suppose that $f,g$ are bounded and $\bE_{f}$ and $\bE_{g}$  are equivalence relations. Suppose $\psi(x) > 0$ $(x > 0)$, and
\ben
\item[$(A_{1})$] there exist $\ve > 0$, $M < \omega$ such that for every $n > M$ and $x,y \in [0,1]$, $$\vp(x) \leq \ve \vp(y) \vp(1/2^n) \Rightarrow x \leq \frac{y}{2^{n+1}};$$ 
\item[$(A_{2})$] $\lim _{n \rar \infty} \psi(1/2^{n})/\vp(1/2^{n}) = 0$.
\een Then $\bE_{g} \not \leq _{B} \bE_{f}$.
\end{theorem}

Observe that $\vp\equiv 1$, $\psi = \Id^{\beta}$ $(0 < \beta < \infty)$ satisfy the assumptions of Theorem \ref{Nr}, so it generalizes \cite[Theorem 2.2 p.\ 1840]{DH}.

The proof  of \cite[Theorem 2.2 p.\ 1840]{DH} has two fundamental constituents. The first idea is to pass to a subspace $X \ss [0,1]^{\omega}$ where a hypothetic Borel reduction $\vt$ of $\bE_{g}$ to $\bE_{f}$ is modular, i.e.\  for $x \in X$, $\vt(x)$ consists of  finite blocks, each of which depends only on a single coordinate of $x$. This technique can be adopted without any difficulty. The second tool is an excessive use of the fact that for $f = \Id^{p}$, $f^{-1}(\|\cdot\|_{f})$ is a norm, which does not follow from the assumptions of Theorem \ref{Nr}. We  get around this difficulty by exploiting that $\vp$ is a perturbation when compared to $\Id^{\alpha}$.

\medskip

\textbf{Proof of Theorem \ref{Nr}.} Suppose $\vt \colon [0,1]^{\omega} \rar [0,1]^{\omega}$ is a Borel reduction of $\bE_{g}$ to $ \bE_{f}$. With $Z_{k} = \{i/2^{k}\colon 0 \leq i \leq 2^{k}\}$, set $Z = \prod_{k < \omega} Z_{k}$; then $\vt$ is a Borel reduction of $\bE_{g}|_{Z \times Z}$ to $\bE_{f}$. For every finite sequence $t \in \prod_{i < |t|} Z_{i}$, let $N_{t} = \{z \in Z \colon z(i) = t(i)~(i < |t|)\}$. We import several lemmas from \cite{DH}.

\begin{lemma}\label{DH1}{\rm (\cite[Claim (i) p.\ 1840]{DH})} For any $j,  k  < \omega$ there exist $l < \omega$, a  finite sequence $s^{\star} \in \prod_{i < |s^{\star}|} Z_{k+i}$, and a comeager set $D \ss  Z$ such that for all $x, \hat x \in D$, if we have $x = r^{\frown}s^{\star \frown} y$
and $\hat x = \hat r ^{\frown} s^{\star \frown} y$ for some $r, \hat r \in  [0,1]^{k}$ and $y \in  [0,1]^{\omega}$, then $\|(\vt (x)- \vt(\hat x))|_{\omega \sm l}\|_{f} < 2^{-j}.$
\end{lemma}
\prf For every $l < \omega$, we define $F_{l} \colon Z \rar \real$ by  \Kem\notag F_{l}(x) = \max\{ \|(\vt (z)- \vt(\hat z))|_{\omega \sm l}\|_{f} \colon  \\ z,\hat z \in Z,~z(i) = \hat z(i) = x(i)~(k \leq i < \omega)\}.\end{multline}
For fixed $x \in Z$, there are only  finitely many $z, \hat z \in Z$ satisfying $z(i) = \hat z(i) = x(i)$ $(k \leq i < \omega)$. For each such pair we have $\|z - \hat z\|_{g} < \infty$, hence  $\|\vt(z) - \vt( \hat z)\|_{f} < \infty$, in particular  $\lim_{l \rar \infty} \|(\vt(z) -  \vt(\hat z))|_{\omega \sm l}\|_{f} = 0$. So $F_{l}(x) < \infty$ for all $l < \omega$ and $\lim_{l \rar \infty}  F_{l}(x) = 0$ $(x \in Z)$. Therefore, by the Baire Category
Theorem, there exists an $l < \omega$ such that $\{x \in Z \colon  F_{l}(x) < 2^{-j}\}$ is not meager. By $f$ being Borel, this set has
the property of Baire, so there is a nonempty open set $O$ on which it is relatively
comeager.

We can assume $O = N_{t}$ for some  finite sequence $t \in \prod_{i < |t|} Z_{i}$, and we can also 
assume $|t| \geq k$. Let $t = r^{\star \frown} s^{\star}$  where $|r^{\star}| = k$. But $F_{l}(x)$ does not depend
on the  first $k$ coordinates of $x$, so $\{x \in Z \colon F_{l}(x) < 2^{-j}\}$ is also relatively comeager
in $N_{r^{\frown}s^{\star}}$  for all $r \in \prod_{i < k } Z_{i}$. Let $D$ be a comeager set such that $F_{l}(x) < 2^{-j}$
whenever $x \in  D  \cap  N_{r^{\frown}s^{\star}}$  for any $r$ of length $k$. Now the conclusion of the claim
follows from the definition of $F_{l}$.$\bs$

\medskip

By \cite[(8.38) Theorem p.\ 52]{K} there is a dense $G_{\delta}$  set $C \ss Z$ such that $\vt|_{C}$  is continuous.

\begin{lemma}\label{DH2}{\rm (\cite[Claim (ii) p.\ 1841]{DH})} For any $j,  k,l   < \omega$ there is a  finite sequence $s^{\star \star} \in \prod_{i < |s^{\star \star}|} Z_{k+i}$ such that for all $x, \hat x \in C$, if we have $x = r^{\frown}s^{\star \star \frown} y$
and $\hat x = r^{\frown} s^{\star\star \frown} \hat y$ for some $r \in  [0,1]^{k}$ and $y , \hat y \in  [0,1]^{\omega}$, then $\|(\vt (x)- \vt(\hat x))|_{l}\|_{f} < 2^{-j}.$

Furthermore, if $G$ is a given dense open subset of $Z$, then $s^{   \star \star}$ can be
chosen such that $N_{r ^{\frown}s^{\star \star}} \ss G$ for all $r \in \prod_{i < k} Z_{i}$.
\end{lemma}
\prf There are only  finitely many $r \in\prod_{i < k } Z_{i}$; enumerate them as $r_{0},  r_{1}, \dots, r_{M-1}$. We construct $s^{\star \star}$ by successive extensions.

Let $t_{0} = \es$. Let $m < M$ and suppose that we have the finite
sequence $t_{m} \in \prod_{i < |t_{m}|} Z_{k+i}$. The basic open
set $N_{r_{m}^{\frown}t_{m}}$ meets the comeager set $C$, so we can pick $w \in C \cap N_{r_{m}^{\frown}t_{m}}$. Since $\vt$ is continuous on $C$ and $f$ is continuous, we can  pass to a smaller open neighborhood $O$ of $w$ such
that for all $x, \hat x \in  C \cap O$, $\|(\vt(x) -  \vt(\hat x))|_{l}\|_{f} < 2^{-j}$. We can assume $O = N_{r_{m}^{\frown}t_{m}'}$ for some extension $t_{m}'$ of $t_{m}$. Since $G$ is dense open, we can further extend $t_{m}'$  to get
$t_{m+1}$ such that $N_{r_{m}^{\frown}t_{m+1}}  \ss  G$.
Once the sequences $t_{m}$ $(m \leq M)$ are constructed, $s ^{\star \star} = t_{M}$ fulfills the requirements.$\bs$

\begin{lemma}\label{DH3}{\rm \cite[Claim (iii) p.\ 1842]{DH}} There exist strictly increasing sequences $(b_{i})_{i < \omega}, (l_{i})_{i < \omega} \ss \omega$ and functions $f_{i} \colon Z_{b_{i}} \rar [0,1]^{l_{j+1}-l_{j}}$ such that $b_{0} = l_{0} = 0$, for $Z' = \prod_{i < \omega} Z_{b_{i}}$ and $\vt' \colon Z' \rar [0,1]^{\omega}$, $\vt'(x) = f_{0}(x_{0}) ^{\frown} \dots ^{\frown}   f_{i}(x_{i})  ^{\frown} \dots$ we have \Keq\label{DH3_5}\|x-\hat x\|_{g} < \infty \Leftrightarrow \|\vt'(x)-\vt'(\hat x)\|_{f} < \infty.\Zeq
\end{lemma}
\prf We construct the sequences $(b_{i})_{i < \omega}, (l_{i})_{i < \omega} \ss \omega$, finite sequences $s_{i} $ $(i < \omega)$ and dense open sets $D^{j}_{i}$ $(i,j < \omega)$ by induction, as follows.

 We have $b_{0} = l_{0} = 0$. Let $j < \omega$ and suppose that we have $b_{j}$, $l_{j}$ and $D_{i}^{j'}$ for every $i < \omega$ and $j' < j$. We apply Lemma \ref{DH1} for $j$ and $k = b_{j}+1$ to get $l_{j+1} = l < \omega$, a finite sequence $s_{j}^{\star} \in \prod_{i < |s_{j}^{\star}|} Z_{b_{j}+1+i}$ and a comeager set $D^{j} \ss Z$ satisfying the conclusions of Lemma \ref{DH1}. We can assume $l_{j+1} > l_{j}$ and $D^{j} \ss C$. Let $(D^{j}_{i})_{i < \omega}$ be a decreasing sequence of dense open subsets of $Z$ such that $ \bigcap_{i < \omega}D^{j}_{i} \ss D^{j}$. We apply Lemma \ref{DH2} for $j$, $k = b_{j} + 1 + |s_{j}^{\star}|$, $l = l_{j+1}$,  and $G = \bigcap_{j' < j} D^{j'}_{j}$ to get $s_{j}^{\star \star}$ as in Lemma \ref{DH2}. We set $s_{j} = s_{j}^{ \star \frown} s_{j}^{\star \star}$ and $b_{j+1} = b_{j} + 1 + |s_{j}|$. 

Let $Z' = \prod_{i < \omega} Z_{b_{i}}$ and set $h \colon Z' \rar Z$, $$h(x) = x_{0} ^{\frown}s_{0}^{\frown}x_{1} ^{\frown}s_{1}^{\frown} \dots ^{\frown}x_{i} ^{\frown}s_{i}^{\frown} \dots . $$ For every $i < \omega$, we define $f_{i} \colon Z_{b_{i}} \rar [0,1]^{l_{j+1}-l_{j}}$ by \Keq\label{FFE}f_{i}(a) = \vt(h(\underbrace{0 ^{\frown} \dots ^{\frown} 0}_{i}~\!\! ^{\frown} a ^{\frown} 0 ^{\frown} 0 ^{\frown} \dots))|_{l_{j+1} \sm l_{j}};\Zeq and we set $\vt' \colon Z' \rar [0,1]^{\omega}$, $\vt'(x) = f_{0}(x_{0}) ^{\frown} \dots ^{\frown}   f_{i}(x_{i})  ^{\frown} \dots$

It remains to prove (\ref{DH3_5}). To see this, it is enough to prove $\|\vt'(x)  - \vt (h(x))\|_{f} < \infty$ for every $x \in Z'$ since then for every $x, \hat x \in Z'$, by $(R_{2})$, \begin{multline} \notag \|\vt'(x) - \vt ' (\hat x)\|_{f}  < \infty \iff \|\vt(h(x)) - \vt  (h(\hat x))\|_{f}  < \infty  \iff \\ \|h(x) -   h(\hat x)\|_{g}  < \infty \iff \|x  - \hat x \|_{g} < \infty .\end{multline}

Let $x \in Z'$ be arbitrary; for every $j < \omega$ we define $e_{j}, e'_{j}  \in Z'$ by setting $$\PR_{i} e_{j} = \b\{ \begin{array}{ll} x_{i}, & \textrm{ if } i = j; \\ 0, & \textrm{ if } i \in \omega \sm \{ j\}; \end{array} \j.,~\PR_{i} e'_{j} = \b\{ \begin{array}{ll} x_{i}, & \textrm{ if } i \leq j; \\ 0, & \textrm{ if }  j < i < \omega. \end{array} \j.$$ Since $h(x)$ and $h(e'_{j})$ agree on all coordinates below $b_{j+1}$, by the definition of $s_{j}^{\star \star}$, $$\|(\vt(h(x)) - \vt(h(e'_{j})))|_{l_{j+1}}\|_{f} < 2^{-j} ~(j < \omega).$$ On the other hand, for $j > 0$, $h(e'_{j})$ and $h(e_{j})$ agree on all coordinates above $b_{j-1}$, so by the definition of $s_{j-1}^{\star}$, \Keq\label{potya}\|(\vt(h(e'_{j})) - \vt(h(e_{j})))|_{\omega \sm l_{j}}\|_{f} < 2^{-j+1} ~(0 <j < \omega).\Zeq Moreover, (\ref{potya}) holds for $j=0$, as well. Then by $(R_{2})$, \begin{multline}\notag \|(\vt'(x)  - \vt (h(x)))|_{l_{j+1} \sm l_{j}}\|_{f}  = \|(\vt(h(e_{j}))  - \vt (h(x)))|_{l_{j+1} \sm l_{j}}\|_{f} \lesssim  \\  \|(\vt(h(e_{j}))  - \vt (h(e'_{j})))|_{\omega\sm l_{j}}\|_{f}+\|(\vt(h(e'_{j}))  - \vt (h(x)))|_{l_{j+1}}\|_{f} \leq 3 \cdot 2^{-j}. \end{multline} Therefore $$\|(\vt'(x)  - \vt (h(x)))\|_{f}  = \sum _{j <\omega} \|(\vt'(x)  - \vt (h(x)))|_{l_{j+1} \sm l_{j}}\|_{f}  \leq \sum _{j <\omega} 3 \cdot 2^{-j}< \infty, $$ as required.$\bs$

\begin{lemma}\label{DH4}{\rm \cite[Claim (iv) p.\ 1843]{DH}} There exist $ c >0$ and $N < \omega$ such that with the notation of (\ref{FFE}), for every $i >N$,
$\| f_{i}(1) - f_{i}(0)\|_{f} > c$.
\end{lemma}
\prf If not, then we can  find a strictly increasing sequence $(j_{m})_{m < \omega} \ss \omega$ such that 
$\|f_{j_{m}}(1) - f_{j_{m}}(0)\|_{f} \leq  2^{-m}$ $(m < \omega)$. Let $\hat x$ be the constant 0 sequence, and let $x$ be the sequence which is 1 at each coordinate $j_{m}$ $(m < \omega)$ and 0 at all other
coordinates. Then $\|x - \hat x\|_{g} = \infty$ but
\begin{multline} \notag \| \vt ' (x) -  \vt ' (\hat x)\|_{f}  = \sum_{j < \omega} \| f_{j}(x(j)) - f_{j}( \hat x(j))\|_{f} = \\ \sum_{m < \omega} \|f_{j_{m}}(1) - f_{j_{m}}(0)\|_{f} \leq \sum_{m < \omega} 2^{-m} < \infty,\end{multline}
contradicting (\ref{DH3_5}).$\bs$

\begin{lemma}\label{DH5} Let $c > 0$, $N < \omega$ be as in Lemma \ref{DH4}. For every $0< D < \omega$ there exists $ N_{D} > \max\{N,D\}$ such that for every $i \geq N_{D}$ there is a $0 \leq k< 2^{b_{N_{D}}}$ with \Keq\label{KE3}\|f_{i}((k+1)/2^{b_{N_{D}}}) - f_{i}(k/2^{b_{N_{D}}}) \|_{f} \geq Dg(1/2^{b_{N_{D}}}).\Zeq
\end{lemma}
\prf Let $\ve > 0$ and $M < \omega$ be as in the assumptions of Theorem \ref{Nr}. Fix $0<D < \omega$; by $(A_{2})$ there exists $N_{D} > \max\{M,N,D\}$ such that with $n = 2^{b_{N_{D}}}$, $2D/c  < \ve \vp(1/n) / \psi(1/n)$. Fix $i \geq N_{D}$, set $l = l_{i+1} - l_{i}$ and $$\gamma_{j} = |\PR_{j}(f_{i}(1) - f_{i}(0))|~(j < l).$$ For every $x= (x_{j})_{j < l} \in [-1,1]^{l}$ set $$\|x\|_{\Delta} = \b(\sum_{j < l} |x_{j}|^{\alpha}\vp(\gamma_{j})\j)^{1/\alpha};$$ then $\|\cdot\|_{\Delta}$ satisfies the triangle inequality on $[0,1]^{l}$. Since $\|f_{i}(1) - f_{i}(0)\|_{f} =\|f_{i}(1) - f_{i}(0)\|_{\Delta} ^{\alpha} = \sum_{j < l}\gamma_{j}^{\alpha}\vp(\gamma_{j})$, by the triangle inequality there is a $0 \leq k < n $ such that $$\|f_{i}((k+1)/n) - f_{i}(k/n) \| _{\Delta}  \geq \frac{1}{n}\|f_{i}(1) - f_{i}(0)\|_{f}^{1/\alpha}.$$ With such a $k$, set $$\delta_{j} = |\PR_{j}(f_{i}((k+1)/n) - f_{i}(k/n))|~ (j < l);$$ i.e.\ we have \Keq\label{KE1}\sum_{j < l} \delta_{j}^{\alpha}\vp(\gamma_{j})\geq \frac{1}{n^{\alpha}}\|f_{i}(1) - f_{i}(0)\|_{f}.\Zeq

Set $J = \{j < l \colon \vp(\gamma_{j}) \leq \vp(\delta_{j})c/(2D\psi(1/n)) \}$. Then \Keq\label{KE2}\sum_{j < l} \delta_{j}^{\alpha}\vp(\gamma_{j}) \leq \sum_{j \in J} \delta_{j}^{\alpha}\vp(\delta_{j})\frac{c}{2D\psi(1/n)} +\sum_{j \notin J} \delta_{j}^{\alpha}\vp(\gamma_{j}).\Zeq By the choice of $N_{D}$, $2D/c  < \ve \vp(1/n) / \psi(1/n)$. So for $j \notin J$, $\vp(\delta_{j}) < \ve\vp(\gamma_{j})\vp(1/n)$. This, by $(A_{1})$ and by $b_{N_{D}} \geq N_{D} >M$, implies $\delta_{j} \leq \gamma_{j}/(2n)$ $(j \notin J)$. Hence $$\sum_{j \notin J} \delta_{j}^{\alpha}\vp(\gamma_{j}) \leq \frac{1}{(2n)^{\alpha}} \sum_{j \notin J} \gamma_{j}^{\alpha}\vp(\gamma_{j})  = 2^{-\alpha}\frac{\|f_{i}(1) - f_{i}(0)\|_{f}}{n^{\alpha}} .$$ So by (\ref{KE1}) and (\ref{KE2}), $$\sum_{j  \in J} \delta_{j}^{\alpha}\vp(\delta_{j})\frac{c}{2D\psi(1/n)} \geq (1-2^{-\alpha})\frac{\|f_{i}(1) - f_{i}(0)\|_{f}}{n^{\alpha}} $$ which implies $$\|f_{i}((k+1)/n) - f_{i}(k/n) \| _{f} = \sum_{j  < l} \delta_{j}^{\alpha}\vp(\delta_{j})\geq D\frac{\psi(1/n)}{n^{\alpha}} = Dg(1/n),$$ as required.$\bs$

\medskip

For every $0<D < \omega$ let $N_{D}$ be as in Lemma \ref{DH5}. Since $g(0) = 0$ and $g$ is continuous, by reassigning $N_{D}$ we can assume $g(1/2^{b_{N_{D}}}) \leq 1/D^{2}$ $(0<D < \omega)$. Let $I_{D} \ss \omega \sm N_{D}$ $(0<D < \omega)$ be pairwise disjoint sets such that $ 1/D \leq |I_{D}| D g(1/2^{b_{N_{D}}}) < 2/D$. For every $0<D < \omega$ and $i \in I_{D}$ pick a $0 \leq k_{i,D}< 2^{b_{N_{D}}}$ satisfying (\ref{KE3}).  Define $x , \hat x  \in Z'$ by $x(i)= k_{i,D}/2^{b_{N_{D}}}$, $\hat x(i)= (k_{i,D}+1)/2^{b_{N_{D}}}$ $(i \in I_{D}, ~0 < D < \omega)$, else $x(i) = \hat x (i) = 0$. Then \Kem\notag \|\vt'(\hat x)- \vt'(x)\|_{f} = \\ \sum_{0 < D < \omega} \sum_{i \in I_{D}} \|f_{i}((k_{i,D}+1)/2^{b_{N_{D}}}) - f_{i}(k_{i,D}/2^{b_{N_{D}}}) \|_{f} \geq \\ \sum_{0 < D < \omega}D|I_{D}|g(1/2^{b_{N_{D}}})= \infty\end{multline} while $$\|\hat x-  x \|_{g} = \sum _{0 < D < \omega} |I_{D}| g(1/2^{b_{N_{D}}}) < \sum_{0 < D < \omega} \frac{2}{D^{2}} < \infty; $$ i.e.\ $x \bE_{g} \hat x$ but $\vt'(x) \not\!\!\bE_{f} \vt'(\hat x)$. This contradiction completes the proof.$\bs$

\section{Applications} \label{con}
 
In this section we construct several families of functions for which our reducibility and nonreducibility results can be applied. Let $1 \leq \alpha \leq \beta< \infty$, let $\vp \colon (0,1] \rar \real$, $\psi \colon [0,1] \rar \real$ be continuous functions and set $f = \Id^{\alpha}\vp$, $f(0) = 0$ and $g = \Id^{\beta}\psi$. 

\subsection{Definition of $\vp$ from $\psi$ and $\kappa$}

In order to facilitate the checking of the conditions of Theorem \ref{fole}, we may use the following approach. 
Instead of defining $\kappa$ from $\vp$ and $\psi$, we may define $\vp$ from $\psi$ and $\kappa$. To this end we set $\kappa(1/2^{n}) = \mu(n)/2^{n\alpha/\beta}$ $(n < \omega)$ where $\mu$ will be specified later. We assume $\mu(0) = \vp(1) = \psi(1)=1$. Then (\ref{U2}), (\ref{U3}) and (\ref{star}) read as \Keq\label{Z1}\vp\b(\frac{1}{2^{n}}\j)  =\sum_{i=0}^{n} 2^{(\alpha-\beta)(n-i)}\mu(i)^{\beta} \psi\b(\frac{2^{(1-\alpha/\beta)i}\mu(i)}{2^{n}}\j) ~(n < \omega),\Zeq \Keq\label{Z2}\sum_{i=n}^{\infty}\frac{1}{2^{i\alpha}} \mu(i)^{\beta}  \psi\b(\frac{\mu(i)}{2^{i\alpha/\beta}}\j) \leq L\sum_{i=0}^{n}2^{i(\beta - \alpha)}\frac{\mu(i)^{\beta} }{2^{n\beta}}  \psi\b(\frac{2^{(1-\alpha/\beta)i}\mu(i)}{2^{n}}\j),\Zeq \Keq\label{Zstar} \mu(n) \leq L\cdot \max_{i < n} \mu(i)\frac{2^{(n-i)\alpha/\beta}}{2^{n-i}} ~(n < \omega).\Zeq Given $\mu$ and $\psi$, we can define $\vp(1/2^{n})$ $(n < \omega)$ by (\ref{Z1}) and then extend $\vp$ to $(0,1]$ to be a continuous function which is affine on $[1/2^{n+1}, 1/2^{n}]$ $(n < \omega)$. When we say below ``we define $\vp$ from $\mu$, $\alpha$, $\beta$ and $\psi$", we mean this definition.  

We show that for a $\vp$ defined this way, if there exist $\ve>0$, $M < \omega$ such that for $n > M$, \Keq\label{Z14}\vp(1/2^{i}) \leq \ve \vp(1/2^{j})\vp(1/2^{n}) \Rightarrow i \geq j+n+3 ~(i,j < \omega)\Zeq then  $(A_{1})$ of Theorem \ref{Nr} holds. Let $x,y \in (0,1]$, say $1/2^{i+1} < x \leq  1/2^{i}$ and $1/2^{j+1} < y \leq 1/2^{j}$. We have $$\vp(x) \in [\vp(1/2^{i+1}), \vp(1/2^{i})],~\vp(y) \in [\vp(1/2^{j+1}), \vp(1/2^{j})],$$ thus $\vp(x) \leq \ve \vp(y)\vp(1/2^{n})$ implies $$ \min\{\vp(1/2^{i+1}), \vp(1/2^{i})\} \leq \ve \max \{\vp(1/2^{j+1}), \vp(1/2^{j})\}\vp(1/2^{n}).$$ So by (\ref{Z14}), for $n > M$ we have $i \geq n+j+2$, which implies $x \leq y/2^{n+1}$, as required. 

\subsection{Explicit examples}

We introduce a family of functions for which our theorems can be applied and whose growth order is easy to calibrate. For $n < \omega$, let $t_{n} \colon (0,1] \rar \real$, $$t_{n}(x) =  \underbrace{1+\log(1+ \dots \log}_{n}(1-\log(x))\dots)~(0 < x \leq 1).$$ For $\eta \in [0,1)^{< \omega}$ we define $l_{\eta} \colon (0,1] \rar \real$, $l_{\eta}(x)= \prod_{i < |\eta|}  t_{i}^{\eta_{i}}~(0 < x \leq 1);$ e.g., \begin{multline}\notag l_{\es}(x) = 1, ~l_{(\eta_{0})}(x) =(1-\log(x))^{\eta_{0}}, \\ l_{(\eta_{0}\eta_{1})}(x) =(1-\log(x))^{\eta_{0}}(1+\log(1-\log(x)))^{\eta_{1}},~\textrm{etc.}\end{multline} Let $<_{\textrm{lex}}$ denote the lexicographic order. We summarize some elementary properties of the functions $l_{\eta}$, which will be used in the sequel. 

\begin{lemma}\label{y} For every $\eta, \eta' \in [0,1)^{< \omega}$ with $\eta <_{\textrm{lex}} \eta'$, $1 \leq \alpha < \infty$ and $\delta > 0$,
\ben[(a)]
\item \label{bEc} $1 \leq  l_{\eta}(xy) \leq l_{\eta}(x)l_{\eta}(y) ~(0 < x,y \leq 1)$;
\item \label{y05} $l_{\eta} \circ \Id^{\delta} \approx  l_{\eta} $ and $ l_{\eta} \lesssim \Id^{-\delta}$;
\item \label{y07} $l_{\eta} (1/2^{n+1}) -l_{\eta} (1/2^{n}) \leq 1$ for every $n < \omega$ sufficiently large;
\item\label{y1} $l_{\eta}$ is continuous and strictly decreasing, moreover if $\eta < _{\mathrm{lex}} \eta'$ then $l_{\eta}/l_{\eta'}$ is strictly  increasing in a neighborhood of 0, so by $l_{\eta}(x)/l_{\eta'}(x) > 0$ $(x > 0)$, $l_{\eta}/l_{\eta'}$ is essentially increasing and $\lim _{x \rar +0} l_{\eta}(x)/l_{\eta'}(x)=0$;
\item\label{y2} $\Id^{\delta}l_{\eta}$ is bounded and $\lim _{x \rar +0} x^{\delta}l_{\eta}(x) = 0$;
\item\label{y3} $f(x)=x^{\delta}l_{\eta}(x)$ $(0 < x < 1)$, $f(0)=0$ is continuous, strictly increasing in a neighborhood of 0, so by $f(x) > 0$ $(x > 0)$, $f$ is essentially increasing; 
\item\label{y4}  $f(x)=x^{\alpha}l_{\eta}(x)$ $(0 < x < 1)$, $f(0)=0$ is continuous, satisfies $(R_{1})$ and $(R_{2})$ hence $\bE_{f}$ is an equivalence relation;
\item\label{y5}  $f(x)=x^{\alpha}/l_{\eta}(x)$  $(0 < x < 1)$, $f(0)=0$ is continuous and strictly  increasing, satisfies $(R_{1})$ and $(R_{2})$ hence $\bE_{f}$ is an equivalence relation;
\item\label{y6} $\vp = 1/ l_{\eta}$ satisfies $(A_{1})$ of Theorem \ref{Nr}.
\een
\end{lemma}
\prf  
It is enough to prove (\ref{bEc}) for $t_{n}$ $(n < \omega)$. We do this by induction on $n$. For $n = 0$, the statement follows from \begin{multline}\notag 1 \leq 1-\log(xy) = 1 - \log(x) - \log(y) \leq  \\ 1 - \log(x) - \log(y) + \log(x)\log(y) = (1 - \log(x))(1 - \log(y)).\end{multline} Let now $n > 1$; then $t_{n} = 1+\log t_{n-1}$, hence $ 1 \leq t_{n}$. By the inductive hypothesis, \begin{multline}\notag t_{n}(xy) = 1+\log t_{n-1}(xy) \leq 1+\log t_{n-1}(x) + \log t_{n-1}(y) \leq  \\ (1+\log t_{n-1}(x)) (1+t_{n-1}(y)) = t_{n}(x)t_{n}(y), \end{multline} as required.

Similarly, it is enough to show  (\ref{y05}) for $t_{n}$ $(n < \omega)$; we use induction on $n$. For $n=0$, the first statement follows from $1-\log(x^{\delta}) = 1-\delta \log(x)$ $(0 < x \leq 1)$, while $ t_{0} \lesssim \Id^{-\delta}$ is elementary analysis. Let now $n > 1$; we have $t_{n} = 1+\log t_{n-1}$. By the inductive hypothesis and $t_{n-1} \geq 1$,  $1+\log (t_{n-1} \circ \Id^{\delta}) \approx 1+\log t_{n-1}$, so the first statement follows. Also by the inductive hypothesis, $1+\log t_{n-1} \lesssim 1-\delta \log \lesssim \Id^{-\delta}$, so the proof is complete.

We show $(l_{\eta} (1/2^{n+1}) -l_{\eta} (1/2^{n}))_{n < \omega}$ is a null sequence; then (\ref{y07}) follows. By elementary analysis, for every $\delta \in [0,1)$ and $m < \omega$, $(t_{m}^{\delta} (1/2^{n+1}) -t_{m}^{\delta}(1/2^{n}))_{n < \omega}$ is a null sequence. Since $l_{\eta}$ is a finite product of $t_{m}^{\delta} $s, the statement follows.

Statements  (\ref{y1}), (\ref{y2}) and (\ref{y3}) are elementary analysis. For (\ref{y4}), $(R_{1})$ is immediate; $(R_{2}a)$ follows from $(x+y)^{\alpha} \lesssim x^{\alpha} +y^{\alpha}$ $(0 \leq x,y \leq 1)$ and $l_{\eta}$ being decreasing; while $(R_{2}b)$ follows from $\Id^{\alpha}l_{\eta}$ being essentially increasing.

Consider now (\ref{y5}). Since $l_{\eta}$ is strictly decreasing, $\Id^{\alpha}/l_{\eta}$ is strictly  increasing. So $(R_{1})$ is immediate and $(R_{2}b)$ holds. To see $(R_{2}a)$, observe that by (\ref{bEc}), for $0<v/2 \leq u \leq v \leq 1$ we have $$l_{\eta}(u) \leq l_{\eta}(v/2) \leq l_{\eta}(1/2) l_{\eta}(v).$$  So for $0 < x , y \leq 1$, $$(x+y)^{\alpha}/ l_{\eta}(x+y) \lesssim l_{\eta}(1/2) (x^{\alpha} / l_{\eta}(x) + y^{\alpha} / l_{\eta}(y)),$$ as required.

It remains to prove (\ref{y6}).  It is enough to show that for every $n < \omega$, $$l_{\eta}(x) \geq l_{\eta}(1/2) l_{\eta}(y) l_{\eta}(1/2^{n}) \Rightarrow x \leq \frac{y}{2^{n+1}}~(i,j < \omega).$$  By (\ref{bEc}), $l_{\eta}(y/2^{n+1}) \leq l_{\eta}(1/2) l_{\eta}(y) l_{\eta}(1/2^{n})$, so since $l_{\eta}$ is decreasing, the statement follows.$\bs$

\begin{corollary}\label{ET} Let $1 \leq \alpha  < \infty$ and let $\eta, \eta' \in  [0,1)^{< \omega}$ satisfy $\eta <_{\textrm{lex}} \eta'$.
\ben
\item\label{ET1} The functions $\psi = l_{\eta}$, $g(x) = x^{\alpha}l_{\eta}(x)$ $(0 < x \leq 1)$, $g(0) = 0$  satisfy the conditions of Theorem \ref{ala}.
\item\label{ET3} The functions $\vp(x) = 1/l_{\eta}(x)$, $\psi(x) = 1/l_{\eta'}(x)$ $(0 < x \leq 1)$, $\vp(0) = \psi(0)=0$ and  $f = \Id^{\alpha} / l_{\eta}$, $g = \Id^{\alpha} / l_{\eta'}$ satisfy the conditions of Theorem \ref{Nr}. 
\een
\end{corollary}
\prf Statement \ref{ET1} follows from (\ref{y1}), (\ref{y2}) and (\ref{y4}) of Lemma \ref{y}. For \ref{ET3},  $\bE_{f}$ and $\bE_{g}$ are equivalence relations by (\ref{y5}) of Lemma \ref{y}; while $(A_{1})$ and $(A_{2})$ follow from  (\ref{y6}) and (\ref{y1}) of Lemma \ref{y}. This completes the proof.$\bs$

\subsection{The counterintuitive case}

In this section we present an example illustrating that the comparison of the growth order of functions does not decide Borel reducibility. Let $\alpha=\beta$ and $\psi \equiv 1$. Then (\ref{Z1}) turns to $\vp(1/2^{n})  =\sum_{i=0}^{n}\mu(i)^{\alpha}$, i.e.\ \Keq\label{Mumu}\mu(n)^{\alpha} = \vp(1/2^{n}) - \vp(1/2^{n-1})~(0 < n < \omega);\Zeq (\ref{Z2}) reads as \Keq\label{G1}\sum_{i=0}^{\infty}\frac{1}{2^{i\alpha}} \mu(n+i)^{\alpha} \leq L \vp(1/2^{n});\Zeq and (\ref{Zstar}) means \Keq\label{Mustar} \mu(n) \leq L \cdot \max_{i < n} \mu(i).\Zeq Since $\mu(n)^{\alpha} \leq \vp(1/2^{n})$, (\ref{G1}) holds if \Keq\label{G2}\sum_{i=0}^{\infty}1/2^{i\alpha} \vp(1/2^{n+i}) \leq L \vp(1/2^{n}).\Zeq

\begin{corollary}\label{HGF}
Let $\vp \colon (0,1] \rar (0,+\infty)$ be an essentially decreasing continuous function such that $\Id^{\alpha}\vp$ is essentially increasing, $ \bE_{\mathrm{Id}^{\alpha}\vp}$ is an equivalence relation, for every $\delta > 0$, $\liminf_{x \rar +0} x^{\delta} \vp(x) = 0$ and (\ref{G2}) holds. Define $\mu$ by (\ref{Mumu}) and suppose  (\ref{Mustar}) holds. Then $\bE_{\mathrm{Id}^{\alpha}}$ and $\bE_{\mathrm{Id}^{\alpha}\vp}$ are Borel equivalent.
\end{corollary}
\prf By Theorem \ref{ala}, $\bE_{\mathrm{Id}^{\alpha}} \leq _{B} \bE_{\mathrm{Id}^{\alpha}\vp}$. By Lemma \ref{EINC}, we can assume in addition that $\vp$ is decreasing. Then the definition of $\mu$ in (\ref{Mumu}) is valid. So by Theorem \ref{fole}, $\bE_{\mathrm{Id}^{\alpha}\vp} \leq _{B} \bE_{\mathrm{Id}^{\alpha}}$.$\bs$  

\medskip

We show that  (\ref{G2}) holds if for some $\ve > 0$, $\Id^{\alpha-\ve} \vp$ is essentially increasing. Then  $$1/2^{(n+i)(\alpha - \ve)} \vp(1/2^{n+i}) \lesssim 1/2^{n(\alpha - \ve)} \vp(1/2^{n})~ (i < \omega),$$ i.e.\ $1/2^{i\alpha} \vp(1/2^{n+i}) \lesssim 1/2^{i\ve} \vp(1/2^{n})$ $(i < \omega)$, so the statement follows.  In particular, by Corollary \ref{ET}.\ref{ET1} and by (\ref{y1}), (\ref{y3}) and (\ref{y4}) of Lemma \ref{y}, $\vp = l_{\eta}$ fulfills these requirements  for every $\eta \in [0,1)^{< \omega}$. By (\ref{y07}) of Lemma \ref{y} and by $\mu(0)=1$, (\ref{Mustar}) also holds for $\vp = l_{\eta}$ $(\eta \in [0,1)^{< \omega})$. That is, $\bE_{\mathrm{Id}^{\alpha}l_{\eta}}$ and $\bE_{\mathrm{Id}^{\alpha}}$ are Borel equivalent. We will see below in (\ref{totref}) that for every $\eta \in [0,1)^{< \omega}$, $\bE_{\mathrm{Id}^{\alpha}} <_{B} \bE_{\mathrm{Id}^{\alpha}/l_{\eta}}$. So the comparison of the growth order of functions does not decide Borel reducibility.

\subsection{The $\alpha < \beta$ case}

Since the previous and following sections contain the analysis of the reducibility of $\bE_{\mathrm{Id}^{\beta}}$ to $\bE_{\mathrm{Id}^{\beta}\psi}$, in the $\alpha < \beta$ case we assume $\psi \equiv 1$.  Then (\ref{Z1}) and (\ref{Z2}) turn to \Keq\label{G3}\vp\b(\frac{1}{2^{n}}\j)  =\sum_{i=0}^{n} 2^{(\alpha-\beta)(n-i)}\mu(i)^{\beta} ~(n < \omega),\Zeq \Keq\label{G4}\sum_{i=0}^{\infty}\frac{1}{2^{i\alpha}} \mu(n+i)^{\beta}   \leq L\sum_{i=0}^{n}2^{(\alpha-\beta)(n-i)}\mu(i)^{\beta}  .\Zeq To satisfy (\ref{G3}), we have to define \Keq \label{MUMM} \mu(n)^{\beta} = \vp(1/2^{n}) - \vp(1/2^{n-1})/2^{\beta-\alpha} ~(0 < n < \omega),\Zeq and then (\ref{G4}) follows from (\ref{G2}). 

\begin{corollary}\label{megm} Let $1 \leq \alpha < \beta < \infty$. Suppose $\vp \colon [0,1] \rar \real^{+}$ is continuous, essentially increasing, $\vp/\Id^{\beta-\alpha}$ is essentially decreasing, $ \bE_{\mathrm{Id}^{\alpha}\vp}$ is an equivalence relation  and for the $\mu$ defined by (\ref{MUMM}), (\ref{Zstar}) holds. Then $\bE_{\Id^{\alpha}\vp} \leq _{B} \bE_{\Id^{\beta}}$.
\end{corollary}
\prf By Lemma \ref{EINC}, we can assume $\vp/\Id^{\beta-\alpha}$ is decreasing, so that (\ref{MUMM}) is valid; while $\vp$ being essentially increasing implies (\ref{G2}). So $\bE_{\Id^{\alpha}\vp} \leq _{B} \bE_{\Id^{\beta}}$ follows from Theorem \ref{fole}.$\bs$

\medskip

The assumptions of Corollary \ref{megm} are affordable: 
\bit
\item[-] if $\vp$ is essentially decreasing, Corollary \ref{HGF} gives the Borel equivalence of $\bE_{\Id^{\alpha}\vp}$ and $\bE_{\Id^{\beta}}$ under suitable assumptions; 
\item[-] in order to not to be in the counterintuitive case, we may assume that  $\vp/\Id^{\beta-\alpha-\delta}$ is decreasing for some $\delta > 0$, so by Corollary \ref{megm},  $\bE_{\Id^{\alpha}\vp} \leq _{B} \bE_{\Id^{\beta-\delta}} < _{B} \bE_{\Id^{\beta}}$;
\eit 
So Corollary \ref{megm} indicates that in the $\alpha < \beta$ case growth order decides Borel reducibility. Moreover, in the next section we will see that in order to guarantee $\bE_{\Id^{\alpha}\vp} \leq _{B} \bE_{\Id^{\beta}}$ by growth order estimates, we need $\Id^{\beta}/(\Id^{\alpha} \vp)$ to be  bounded; the assumptions of Corollary \ref{megm} reflect this constraint.

Finally we check that for every $\eta \in [0,1)^{< \omega}$, the function $\vp(0) = 0$, $\vp(x) = 1/l_{\eta}(x)$ $(0 < x \leq 1)$ satisfies the assumptions of Corollary \ref{megm}. By 
$\lim_{n \rar \infty} l_{\eta}(1/2^{n+1})/l_{\eta}(1/2^{n})=1$ we have $$\mu(n+1) \lesssim 1/l_{\eta}^{1/\beta}(1/2^{n+1}) \lesssim 1/l_{\eta}^{1/\beta}(1/2^{n})\lesssim \mu(n) ~(n < \omega),$$ i.e.\ (\ref{Zstar}) holds. The other assumptions follow from (\ref{y1}), (\ref{y3}) and (\ref{y5}) of Lemma \ref{y}. So \Keq\label{totref}\bE_{\Id^{\alpha}/l_{\eta}} \leq_{B} \bE_{\Id^{\gamma}} < _{B} \bE_{\Id^{\beta}}~ (\eta \in [0,1)^{< \omega},~ 1 \leq \alpha < \gamma < \beta < \infty).\Zeq

\subsection{The $\alpha = \beta$ case}

This is the most interesting case for us. Now (\ref{Z1}), (\ref{Z2}) and (\ref{Zstar}) turn to \Keq\label{G5}\vp\b(\frac{1}{2^{n}}\j)  =\sum_{i=0}^{n} \mu(i)^{\alpha} \psi\b(\frac{\mu(i)}{2^{n}}\j) ~(n < \omega),\Zeq \Keq\label{G6}\sum_{i=0}^{\infty}\frac{1}{2^{i\alpha}} \mu(n+i)^{\alpha}  \psi\b(\frac{\mu(n+i)}{2^{n+i}}\j) \leq L\sum_{i=0}^{n}\mu(i)^{\alpha}  \psi\b(\frac{\mu(i)}{2^{n}}\j) ~(n < \omega),\Zeq \Keq\label{Gstar} \mu(n) \leq L \cdot \max_{i < n} \mu(i) ~(n < \omega).\Zeq
 We obtain a sufficient condition for (\ref{G6}) and (\ref{Gstar}).

\begin{lemma}\label{2eset} Assume $\psi$ is essentially increasing, $\psi(x) > 0$ for $x > 0$ and $\mu(n) \leq 1$ for every $n < \omega$ sufficiently large. Then (\ref{G6}) and (\ref{Gstar}) hold.
\end{lemma}
\prf Since $\psi$ is essentially increasing, for every $n$ sufficiently large we have $\psi(\mu(n+i)/2^{n+i})  \lesssim \psi(1/2^{n}) $ $(0 \leq i < \omega)$. Hence   \Keq\label{Z23}\sum_{i=0}^{\infty}\frac{1}{2^{i\alpha}} \mu(n+i)^{\alpha}  \psi\b(\frac{\mu(n+i)}{2^{n+i}}\j) \lesssim \frac{1}{(1-1/2^{\alpha})} \psi(1/2^{n})~(n < \omega),\Zeq thus by $\mu(0) = 1$, (\ref{G6}) follows. Also by $\mu(0)=1$ we have (\ref{Gstar}), so the proof is complete.$\bs$

\subsubsection{The question of S.\ Gao}

In this section, in the spirit of (\ref{nov}), we give the negative answer to the question of S.\ Gao mentioned in the introduction. 

\begin{corollary}\label{Rc} Let $1 \leq \alpha < \infty$ be arbitrary. Let $\mu \colon \omega \rar [0, \infty)$ be such that $\mu(0)=1$. Let $\psi \colon [0,1] \rar [0,\infty)$ be a continuous essentially increasing function such that $\psi(1)=1$, (\ref{G6}) holds and there is a $K > 0$ for which \Keq\label{Z12}\frac{1}{K}\psi(1/2^{n}) \leq \psi\b(\frac{\mu(i)}{2^{n}}\j) \leq K \psi(1/2^{n})~( 0 \leq i \leq n < \omega).\Zeq

Set $\sigma_{\mu^{\alpha}}(n) = \sum_{i=1}^{n}\mu^{\alpha}(i)$ $(n < \omega)$. Suppose $(\b(1+ \sigma_{\mu^{\alpha}}(n) \j) \psi\b(1/2^{n}\j))_{n < \omega}$ is essentially decreasing and (\ref{Gstar}) holds.

Define $\vp$ from $\mu$, $\alpha$  and $\psi$.  Set $f(x) = x^{\alpha}\vp(x)$ $(0 < x \leq 1)$, $f(0)=0$ and $g= \Id^{\alpha}\psi$ and suppose $\bE_{f}$ and $\bE_{g}$  are equivalence relations. Then $\bE_{f} \leq _{B} \bE_{g}$.

If, in addition, $\vp$ satisfies $A_{1}$ of Theorem \ref{Nr}, or equivalently $\vp$ satisfies (\ref{Z14}), and $\lim_{n \rar \infty} \sigma_{\mu^{\alpha}}(n) = \infty$, then $\bE_{g} \not \leq _{B} \bE_{f}$.
\end{corollary}
\prf By (\ref{Z12}), from (\ref{G5}) we get \Keq\label{Z113}\frac{1}{K}\b(1+ \sigma_{\mu^{\alpha}}(n) \j) \psi\b(\frac{1}{2^{n}}\j) \leq \vp\b(\frac{1}{2^{n}}\j)  \leq K\b(1+ \sigma_{\mu^{\alpha}}(n) \j) \psi\b(\frac{1}{2^{n}}\j).\Zeq Since $((1+ \sigma_{\mu^{\alpha}}(n)) \psi(1/2^{n}))_{n < \omega}$ is essentially decreasing, $\vp$ is essentially increasing. So by Theorem \ref{fole}, $\bE_{f} \leq _{B} \bE_{g}$.

Moreover, if $\vp$ satisfies $A_{1}$ of Theorem \ref{Nr}, which follows e.g.\ if $\vp$ satisfies (\ref{Z14}), then since $\lim_{n \rar \infty} \sigma_{\mu^{\alpha}}(n) = \infty$ implies $(A_{2})$ of Theorem \ref{Nr}, $\bE_{g} \not \leq _{B} \bE_{f}$. This completes the proof.$\bs$

\medskip

 Many natural functions satisfy the conditions of Corollary \ref{Rc} for both $\vp$ and $\psi$, in particular the functions $1/l_{\eta}$. By Lemma \ref{mege}, the following result gives the negative answer to the question of S.\ Gao.
 
 \begin{corollary}\label{logos}  For every $1 \leq \alpha < \beta < \infty$ and $\eta, \eta' \in [0,1)^{< \omega}$ with $\eta <_{\textrm{lex}} \eta'$, $$\bE_{\Id^{\alpha}} < _{B} \bE_{\Id^{\alpha}/l_{\eta}} < _{B} \bE_{\Id^{\alpha}/l_{\eta'}}  <_{B} \bE_{\Id^{\beta}}.$$
 \end{corollary}
 \prf By Lemma \ref{y} (\ref{y5}), $\bE_{\Id^{\alpha}/l_{\eta}}$ is an equivalence relation, and in (\ref{totref}) we obtained $\bE_{\Id^{\alpha}/l_{\eta}}<_{B} \bE_{\Id^{\beta}}.$ By Lemma \ref{y} (\ref{y1}), $(l_{\eta'}(1/2^{n})/l_{\eta}(1/2^{n}))_{n < \omega}$ is strictly increasing for $n$ sufficiently large. Thus there is a function $\mu \colon \omega \rar \real^{+}$ such that $\mu(0)=1$, and for every $n < \omega$ sufficiently large, 
  $$\mu^{\alpha}(n) = l_{\eta'}(1/2^{n})/l_{\eta}(1/2^{n}) - l_{\eta'}(1/2^{n-1})/l_{\eta}(1/2^{n-1}).$$
Let $\psi=1/l_{\eta'}$ $(0 < x \leq 1)$, $\psi(0)=0$ and define $\vp$ from $\mu$, $\alpha$ and $\psi$. We check the conditions of Corollary \ref{Rc}.
 
 First we show that for every $\ve > 0$, \Keq\label{UTS}2^{-n\ve} \leq \mu^{\alpha}(n) \leq 1 \Zeq holds for $n$ sufficiently large.   By Lemma \ref{y} (\ref{bEc}), (\ref{y07}) and  (\ref{y1}), for every $n$ sufficiently large, \begin{multline}\notag \mu^{\alpha}(n) = \frac{l_{\eta'}(1/2^{n})}{l_{\eta}(1/2^{n})}-\frac{l_{\eta'}(1/2^{n-1})}{l_{\eta}(1/2^{n-1})} = \\  \frac{l_{\eta'}(1/2^{n})-l_{\eta'}(1/2^{n-1})}{l_{\eta}(1/2^{n})} + \frac{l_{\eta'}(1/2^{n-1})}{l_{\eta}(1/2^{n})} - \frac{l_{\eta'}(1/2^{n-1})}{l_{\eta}(1/2^{n-1})} \leq \\ \frac{l_{\eta'}(1/2^{n})-l_{\eta'}(1/2^{n-1})}{l_{\eta}(1/2^{n})} \leq \frac{1}{l_{\eta}(1/2^{n})} \leq 1. \end{multline} For the lower bound, take an $m > |\eta'|$ and consider $t_{m}$. By Lemma \ref{y} (\ref{y1}), $l_{\eta'}(1/2^{n})/(l_{\eta}(1/2^{n})t_{m}(1/2^{n}))$ is still strictly increasing for $n$ sufficiently large. So for $n$ sufficiently large, $$\mu^{\alpha}(n) = \frac{l_{\eta'}(1/2^{n})}{l_{\eta}(1/2^{n})}-\frac{l_{\eta'}(1/2^{n-1})}{l_{\eta}(1/2^{n-1})} \geq \frac{l_{\eta'}(1/2^{n-1})}{l_{\eta}(1/2^{n-1})}\frac{t_{m}(1/2^{n}) - t_{m}(1/2^{n-1})}{t_{m}(1/2^{n-1})}.$$ It is elementary analysis that $t_{m}(1/2^{n})-t_{m}(1/2^{n-1}) \geq 1/n^{2}$ for $n$ sufficiently large, so the statement follows.

 By Lemma \ref{y} (\ref{y1}), $\psi$ is continuous, essentially increasing and $\psi(1) = 1$. Lemma \ref{2eset} gives (\ref{G6}) and (\ref{Gstar}). Also, (\ref{Z12}) follows from Lemma \ref{y} (\ref{y05}) using that $2^{-n/2} \leq \mu^{\alpha}(n) \leq 2^{n/2}$ holds for $n$ sufficiently large.

We have $$\b(1+ \sigma_{\mu^{\alpha}}(n) \j) \psi\b(1/2^{n}\j) \approx 1/l_{\eta}(1/2^{n})~(n < \omega),$$ so $(\b(1+ \sigma_{\mu^{\alpha}}(n) \j) \psi\b(1/2^{n}\j))_{n < \omega}$ is essentially decreasing.
 
 By (\ref{Z113}), $$\vp\b(\frac{1}{2^{n}}\j) \approx  (1+\sigma_{\mu^{\alpha}}(n))\psi\b(\frac{1}{2^{n}}\j) \approx l_{\eta'}(1/2^{n})/l_{\eta}(1/2^{n}) \psi\b(\frac{1}{2^{n}}\j) \approx 1/l_{\eta}(1/2^{n}),$$ so by Corollary \ref{Rc}, $\bE_{\Id^{\alpha}/l_{\eta}} \leq _{B} \bE_{\Id^{\alpha}/l_{\eta'}}$.
 
 By  Lemma \ref{y} (\ref{y1}), $\lim_{x \rar +0} l_{\eta'}(x)/l_{\eta}(x) = \infty$, i.e.\ $\lim_{n \rar \infty} \sigma_{\mu^{\alpha}}(n) = \infty$. By  Lemma \ref{y} (\ref{y6}), $1/l_{\eta}$ satisfies $A_{1}$ of Theorem \ref{Nr}, so again  by Corollary \ref{Rc}, $\bE_{\Id^{\alpha}/l_{\eta'}} \not \leq _{B} \bE_{\Id^{\alpha}/l_{\eta}} $. The $\eta = \es$ special case gives $\bE_{\Id^{\alpha}}< _{B} \bE_{\Id^{\alpha}/l_{\eta}} $, so the proof is complete.$\bs$
 
 \subsubsection{Embedding long linear orders}
 
In this section we show that every linear order which can be embedded into $(\mc{P}(\omega)/\mathrm{fin}, \subset)$ also embeds into the set of Borel equivalence relations $\bE_{f}$ satisfying $\bE_{\Id^{\alpha}} \leq _{B} \bE_{f} \leq _{B} \bE_{\Id^{\alpha}/(1-\log)}$ ordered by $<_{B}$. We refer to \cite{B-U} for results on embedding ordered sets into $(\mc{P}(\omega)/\mathrm{fin}, \subset)$, here we only remark that it is consistent with ZFC, e.g.\ under the Continuum Hypothesis, that every ordered set of size continuum embeds into $(\mc{P}(\omega)/\mathrm{fin}, \subset)$.
 
 \begin{corollary}\label{om} Let $1 \leq \alpha < \infty$ be fixed. There is a mapping $\mc{F} \colon \mc{P}(\omega)/\mathrm{fin} \rar C[0,1]$ such that for every $U, V \in \mc{P}(\omega)/\mathrm{fin}$, $\bE_{\mc{F}(U)}$ is an equivalence relation satisfying $\bE_{\Id^{\alpha}/(1-\log)^{1-\log(17/16)}} \leq _{B} \bE_{\mc{F}(U)} \leq _{B} \bE_{\Id^{\alpha}/(1-\log)}$ and $U \subset V \Rightarrow \bE_{\mc{F}(V)} <_{B}  \bE_{\mc{F}(U)}$.
\end{corollary} 
\prf Let $ \gamma= 17/16$. For every $U \in \mc{P}(\omega)$ set $$\mu_{U}(0) =1, ~\mu^{\alpha}_{U}(n) = \gamma^{|U \cap \lfloor 1+ \log (n)\rfloor |} ~(0 < n < \omega).$$ Let $\psi_{0}(x) = 1/(1-\log(x))^2$ $(0< x \leq 1)$, $\psi_{0}(0) = 0$. For every $U \in  \mc{P}(\omega)$ we define $\vp_{U}$ from $\mu_{U}$, $\alpha$ and $\psi_{0}$, and we set $\mc{F}(U) = \Id^{\alpha} \vp_{U}$. 

First we show that for every $U \in \mc{P}(\omega)$, $\mu_{U}$ and $\psi_{0}$ satisfy (\ref{Z12}). By definition, $1 \leq \mu^{\alpha}_{U}(n) \leq \gamma^{1+\log(n)} \leq \gamma n$ $(0< n < \omega)$, so (\ref{Z12}) follows.

Next we show that for every $U \in \mc{P}(\omega)$, $\vp_{U}$ is essentially increasing. Since (\ref{Z12}) holds, by (\ref{Z113}) it is enough to show that $((1+\sigma_{\mu^{\alpha}_{U}}(n))\psi_{0}(1/2^{n}))_{n < \omega}$ is essentially decreasing. We have $\psi_{0}(1/2^{n}) \approx 1/n^{2}$ $(0 < n < \omega)$. Let $0 < n < m < \omega$ be fixed, say $m = \rho n$ for some $\rho > 1$. If $\mu_{U}^{\alpha}(n) = \gamma^{k}$, then $\sigma_{\mu^{\alpha}_{U}}(n) \geq n\gamma^{k-1}/2$, and $$\sigma_{\mu^{\alpha}_{U}}(m) \leq \sigma_{\mu^{\alpha}_{U}}(n) + (m-n) \gamma^{k+1+\log(m) - \log(n) } \leq \sigma_{\mu^{\alpha}_{U}}(n) + (\rho-1)n \gamma^{k+1+\log(\rho)}.$$ Hence \begin{multline} \notag \frac{\sigma_{\mu^{\alpha}_{U}}(m)}{m^{2}} \leq \frac{\sigma_{\mu^{\alpha}_{U}}(n) + (\rho-1) n \gamma^{k+1+\log(\rho) }}{(\rho n )^{2}} \leq \\ \frac{\sigma_{\mu^{\alpha}_{U}}(n)}{n^{2}} + \frac{\gamma^{k+1+\log(\rho) }}{\rho n } \leq \frac{\sigma_{\mu^{\alpha}_{U}}(n)}{n^{2}} \b(1+ \frac{2 \gamma^{2+\log(\rho) }}{\rho} \j) \leq 9 \frac{\sigma_{\mu^{\alpha}_{U}}(n)}{n^{2}}. \end{multline} This shows $((1+\sigma_{\mu^{\alpha}_{U}}(n))\psi_{0}(1/2^{n}))_{n < \omega}$ is essentially decreasing.

Next we check that  for every $U \in \mc{P}(\omega)$, $\bE_{\mc{F}(U)}$ is an equivalence relation. By definition, $(R_{1})$ holds; $(R_{2}a)$ holds for $\Id^{\alpha}\psi_{0}$ with $C=8\alpha$, so since $\vp_{U}/\psi_{0}$ is decreasing, $(R_{2}a)$ holds for $\Id^{\alpha}\vp_{U}$, as well. Finally $(R_{2}b)$ follows from $\Id^{\alpha}\vp_{U}$ is essentially increasing. 

Our task is to prove that if $U,V \in \mc{P}(\omega)$ satisfy $U \ss ^{\star} V$, $|V \sm U| = \infty$ then $\bE_{\mc{F}(V)} <_{B}  \bE_{\mc{F}(U)}$. Observe that if $U, U' \in \mc{P}(\omega)$ differ only by a finite set then  $\mc{F}(U) \approx \mc{F}(U')$ hence $\bE_{\mc{F}(U)}=\bE_{\mc{F}(U')}$. So we can assume $U \ss V$, $0 \in V \sm U$. 

Our strategy is to  show that $\vp=\vp_{V}$ can be obtained from $\psi = \vp_{U}$ as in (\ref{G5}) with a $\mu$ satisfying the assumptions of Corollary \ref{Rc}. Set $\mu(0)=1$, $$\mu^{\alpha}(n+1) = \frac{1+\sigma_{\mu^{\alpha}_{V}}(n+1)}{1+\sigma_{\mu^{\alpha}_{U}}(n+1)} - \frac{1+\sigma_{\mu^{\alpha}_{V}}(n)}{1+\sigma_{\mu^{\alpha}_{U}}(n)} ~(n < \omega).$$ Later on will prove \Keq\label{F1}\frac{\gamma-1}{(n+2)^3} \leq \mu^{\alpha}(n) \leq 1 ~(n < \omega);\Zeq now we assume (\ref{F1}) and verify the conditions of Corollary \ref{Rc}.

We have $\vp_{U}$ is continuous and $\vp_{U}(0) = 1$. As we have seen above, $\vp_{U}$ is essentially increasing. By $\mu \leq 1$, Lemma \ref{2eset} gives (\ref{G6}) and (\ref{Gstar}). By (\ref{F1}),  $$\vp_{U}\b(\frac{1}{2^{n+3\lfloor \log(n+2)\rfloor +7}}\j) \lesssim \vp_{U}\b(\frac{\mu(i)}{2^{n}}\j) \lesssim \vp_{U}\b(\frac{1}{2^{n}}\j) ~(0 \leq i \leq n < \omega),$$ so (\ref{Z12}) follows from \begin{multline} \notag\vp_{U}\b(1/2^{n}\j) \approx (1+\sigma_{\mu^{\alpha}_{U}}(n))\psi_{0}(1/2^{n})\approx \\ (1+\sigma_{\mu^{\alpha}_{U}}(n+3\lfloor \log(n+2)\rfloor +7))\psi_{0}(1/2^{n+3\lfloor \log(n+2)\rfloor +7}) \approx \\ \vp_{U}\b(1/2^{n+3\lfloor \log(n+2)\rfloor +7}\j).\end{multline}

Let $\vp$ be defined from $\mu$, $\alpha$ and $\vp_{U}$. We have $$1+ \sigma_{\mu^{\alpha}}(n) = (1+\sigma_{\mu^{\alpha}_{V}}(n))/(1+\sigma_{\mu^{\alpha}_{U}}(n))~(n < \omega),$$ so by (\ref{Z113}), \begin{multline} \notag \vp\b(\frac{1}{2^{n}}\j) \approx  \frac{1+\sigma_{\mu^{\alpha}_{V}}(n)}{1+\sigma_{\mu^{\alpha}_{U}}(n)} \vp_{U}(n) \approx \\ \frac{1+\sigma_{\mu^{\alpha}_{V}}(n)}{1+\sigma_{\mu^{\alpha}_{U}}(n)} (1+\sigma_{\mu^{\alpha}_{U}}(n)) \psi_{0}\b(\frac{1}{2^{n}}\j) = \vp_{V}\b(\frac{1}{2^{n}}\j) .\end{multline} Thus $\bE_{\mc{F}(V)}= \bE_{\Id^{\alpha}\vp}$; and $(\b(1+ \sigma_{\mu^{\alpha}}(n) \j) \psi\b(1/2^{n}\j))_{n < \omega}$ is essentially decreasing. So by Corollary \ref{Rc}, $\bE_{\mc{F}(V)}   \leq_{B} \bE_{\mc{F}(U)}$.

Observe that $\psi_{0}$ satisfies (\ref{Z14}) with $M=0$ and $\ve = 1/8$. Since $(1+\sigma_{\mu^{\alpha}_{U}}(n))_{n < \omega}$ is increasing, $\vp_{U}(1/2^{n}) \approx (1+\sigma_{\mu^{\alpha}_{U}}(n))\psi_{0}(1/2^{n})$ $(n < \omega)$ also satisfies (\ref{Z14}) with the same $M$ and a smaller $\ve$. Thus $\vp_{U}$ satisfies $A_{1}$ of Theorem \ref{Nr}.

Since $  |U\cap \log(n)|+k \leq  |V \cap \log(n)|$ implies $\gamma^{k} \mu^{\alpha}_{U}(n) \leq \mu^{\alpha}_{V}(n)$, we have $$\lim_{n \rar \infty} (1+\sigma_{\mu^{\alpha}_{U}}(n))/(1+\sigma_{\mu^{\alpha}_{V}}(n)) = 0$$ hence $\lim_{n \rar \infty} \sigma_{\mu^{\alpha}}(n) = \infty$. So again by Corollary \ref{Rc}, $\bE_{\mc{F}(U)}  \not \leq_{B} \bE_{\mc{F}(V)}$. For $U = \es$ and $V = \omega$, $\vp_{U} \approx 1/(1-\log)$ and $\vp_{V} \approx 1/(1-\log)^{1-\log(\gamma)}$, so $\bE_{\Id^{\alpha}/(1-\log)^{1-\log(\gamma)}} \leq _{B} \bE_{\mc{F}(U)} \leq _{B} \bE_{\Id^{\alpha}/(1-\log)}$ $(U \in  \mc{P}(\omega))$.

It remains to prove (\ref{F1}). For $n=1$, $\mu^{\alpha}(1) = (1+\gamma)/2 - 1$; for $n=2$, $\mu^{\alpha}(2) = (1+2\gamma)/3 - (1+\gamma)/2$. So (\ref{F1}) holds for $n=1,2$. Let $n \geq 2$; then $$1+\sigma_{\mu^{\alpha}_{V}}(n) =1+\gamma + a_{n}, ~1+\sigma_{\mu^{\alpha}_{U}}(n) =2 + b_{n}$$ and $$1+\sigma_{\mu^{\alpha}_{V}}(n+1) =1+\gamma + a_{n}+\gamma^{c}, ~1+\sigma_{\mu^{\alpha}_{U}}(n+1) =2 + b_{n}+\gamma^{d},$$ where $c \geq d+1$ and $\gamma \leq a_{n}/b_{n}\leq \gamma^{c-d}$  $(1 <n  < \omega)$. Then for every $2 \leq n < \omega$, \begin{multline}\label{F2}\frac{1+\sigma_{\mu^{\alpha}_{V}}(n+1)}{1+\sigma_{\mu^{\alpha}_{U}}(n+1)} - \frac{1+\sigma_{\mu^{\alpha}_{V}}(n)}{1+\sigma_{\mu^{\alpha}_{U}}(n)} = \frac{1+\gamma + a_{n}+\gamma^{c}}{2 + b_{n}+\gamma^{d}} - \frac{1+\gamma + a_{n}}{2 + b_{n}} = \\ \frac{\gamma^{c} - \gamma^{d}\frac{1+\gamma+a_{n}}{2+b_{n}}}{2+b_{n}+\gamma^{d}} =  \frac{\gamma^{c} - \gamma^{d}\b(\frac{a_{n}}{b_{n}} - \frac{\frac{2a_{n}}{b_{n}} - (1 + \gamma)}{2+b_{n}}\j)}{2+b_{n}+\gamma^{d}} \geq  \\\frac{\gamma^{c} - \gamma^{d}\b(\gamma^{c-d} - \frac{2\gamma - (1 + \gamma)}{2+b_{n}}\j)}{2+b_{n}+\gamma^{d}} =\gamma^{d} \frac{\gamma-1}{(2+b_{n})(2+b_{n}+\gamma^{d})} .\end{multline} We have $d \leq \lfloor \log(n)\rfloor$, so $b_{n} \leq n\gamma^{d} \leq n^2$ $(2 \leq n < \omega)$. So (\ref{F2}) can be estimated from below by $$\frac{\gamma-1}{(2/\gamma^{d}+b_{n}/\gamma^{d})(2+b_{n}+\gamma^{d})} \geq \frac{\gamma-1}{(2/\gamma^{d}+n)(2+n^{2}+n)} \geq \frac{\gamma-1}{(n+2)^{3}},$$ as stated. 

For the upper bound, as we have seen in (\ref{F2}), it is enough to show $$\gamma^{c} - \gamma^{d}\frac{1+\gamma+a_{n}}{2+b_{n}} \leq 2+b_{n}+\gamma^{d}~(2 \leq n < \omega).$$ Since $c \leq 1+\log(n+1)$ and $n-2 \leq b_{n}$, $$\gamma^{c} \leq \gamma (n+1)^{\log(\gamma)} \leq n \leq 2+b_{n}~(2 \leq n < \omega)$$ for our $\gamma$. This completes the proof.$\bs$

\end{document}